\documentclass[11pt]{article}
\usepackage{tikz}
\usetikzlibrary{arrows}
\usepackage[T1]{fontenc}
\usepackage{textcomp}
\usepackage{color}

\usepackage[scaled=0.92]{helvet}
\usepackage{times}
\usepackage{amssymb}
\usepackage[all,knot,arc,curve,color,frame]{xy}
\usepackage{url}
\input xy
\xyoption{all}
\usepackage{extarrows}
\usepackage{enumitem}
\setlist{nosep}

\usepackage[all,knot,arc]{xy}
\usepackage{mathtools}

\newtheorem{theorem}{Theorem}[section]
\newtheorem{prop}[theorem]{Proposition}
\newtheorem{lemma}[theorem]{Lemma}
\newtheorem{cor}[theorem]{Corollary}
\newtheorem{conjecture}[theorem]{Conjecture}
\newtheorem{unnumber}{}

\newtheorem{prepf}[unnumber]{Proof}
\newenvironment{pf}{\prepf\rm}{\endprepf}
\newtheorem{preprob}[theorem]{Problem}
\newenvironment{prob}{\preprob\rm}{\endpreprob}

\newcommand{\qed}{\hfill$\Box$}

\newcommand{\pgl}{\mathop{\mathrm{PGL}}}
\renewcommand{\wr}{\mathbin{\mathrm{wr}}}
\newcommand{\cart}{\mathbin{\square}}

\DeclareMathOperator\gr{Gr}
\DeclareMathOperator\Gr{Gr}
\newcommand\Gam{\Gamma}
\newcommand{\End}{\mathop{\mathrm{End}}}
\newcommand{\Aut}{\mathop{\mathrm{Aut}}}

\newcommand\gap{\textsf{GAP}}
\newcommand\magma{{\sc Magma}}

\begin{document}
\title{Primitive groups, graph endomorphisms and synchronization}
\author{Jo\~ao Ara\'ujo\\
  {\small Universidade Aberta, R. Escola Polit\'{e}cnica, 147}\\
  {\small 1269-001 Lisboa, Portugal}\\{\footnotesize \&}\\
  {\small CAUL/CEMAT, Universidade de Lisboa}\\
  {\small 1649-003 Lisboa, Portugal}\\
  {\small jaraujo@ptmat.fc.ul.pt}\\\\
  Wolfram Bentz\\
  {\small CAUL/CEMAT, Universidade de Lisboa}\\
  {\small 1649-003 Lisboa, Portugal}\\
  {\small wfbentz@fc.ul.pt}\\\\
  Peter J. Cameron\\
  {\small Mathematical Institute, University of St Andrews}\\
  {\small North Haugh, St Andrews KY16 9SS, UK}\\
  {\small  pjc20@st-andrews.ac.uk}\\\\
  Gordon Royle\\
  {\small Centre for the Mathematics of Symmetry and Computation}\\
  {\small The University of Western Australia}\\
  {\small Crawley, WA 6009, Australia}\\
  {\small gordon.royle@uwa.edu.au}\\\\
  Artur Schaefer\\
  {\small Mathematical Institute, University of St Andrews}\\
  {\small North Haugh, St Andrews KY16 9SS, UK}\\
  {\small  as305@st-andrews.ac.uk}}
\date{}

\maketitle
\begin{abstract}
Let $\Omega$ be a set of cardinality $n$, $G$ a permutation group on $\Omega$,
and $f:\Omega\to\Omega$ a map which is not a permutation. We say that $G$
\emph{synchronizes} $f$ if the transformation semigroup $\langle G,f\rangle$ contains a constant
map, and that $G$ is a \emph{synchronizing group} if $G$ synchronizes \emph{every} non-permutation.

A synchronizing group is necessarily primitive, but there are primitive groups that
are not synchronizing. Every non-synchronizing primitive group fails to synchronize
at least one uniform transformation (that is, transformation whose kernel has parts of
equal size), and it had previously been conjectured that this was essentially
the only way in which a primitive group could fail to be synchronizing --- in other
words, that a primitive group synchronizes every non-uniform transformation.

The first goal of this paper is to prove that this conjecture is false, by
exhibiting primitive groups that fail to synchronize specific non-uniform transformations of
ranks $5$ and $6$. As it has previously been shown that primitive groups
synchronize every non-uniform transformation of rank at most $4$, these examples
are of the lowest possible rank. In addition we produce graphs with primitive automorphism groups that have approximately $\sqrt{n}$ \emph{non-synchronizing ranks}, thus refuting another conjecture on the number of non-synchronizing ranks of a primitive group.

The second goal of this paper is to extend the spectrum of ranks for which it is known that primitive groups
synchronize every non-uniform transformation of that rank. It has previously been shown that a
primitive group of degree $n$ synchronizes every non-uniform transformation of rank $n-1$ and $n-2$, and here
this is extended to $n-3$ and $n-4$.

Determining the exact spectrum of ranks for which there exist non-uniform transformations not synchronized
by some primitive group is just one of several  natural, but possibly difficult, problems on automata, primitive groups,
graphs and computational algebra arising from this work; these are outlined in the final section.
\end{abstract}

\section{Introduction}

Let $\Omega$ be a set of size $n$ and let $f$ be a transformation on $\Omega$ of rank (size of image) smaller than $n$ (in other words,
$f$ is a non-permutation). A permutation group $G$ of degree $n$ on $\Omega$ \emph{synchronizes}
$f$ if the transformation semigroup $\langle G, f\rangle$ contains a constant transformation. The \emph{kernel} of $f$ is
the partition of $\Omega$ determined by the equivalence relation $x \equiv y$ if and only if $xf = yf$. If
the parts of the kernel all have the same size, then $f$ is called \emph{uniform}; it is
\emph{non-uniform} otherwise. A group is called \emph{synchronizing} if it synchronizes every non-permutation.
A synchronizing group is necessarily primitive (see \cite{ArnoldSteinberg}) but the converse
does not hold and considerable efforts have been made to determine exactly which primitive groups are synchronizing.

To show that a primitive group is \emph{not synchronizing} it is necessary to find a
\emph{witness}, which is a transformation $f$ such that $\langle G, f \rangle$ does
\emph{not} contain a constant map. Neumann \cite{neumann} proved that any non-synchronizing primitive group
has a \emph{uniform} witness. This prompted the definition of a primitive group as
being \emph{almost synchronizing} if it synchronizes every non-uniform
transformation. In \cite{abc} the \emph{almost synchronizing conjecture} was stated. 
This asserts that primitive groups have no non-uniform witnesses; 
that is, they are almost synchronizing.

This conjecture has previously been proved for transformations of
ranks $2$, $3$, $4$, $n-2$ and $n-1$ or, in other words, for
transformations of \emph{very low}, or \emph{very high}, rank (see \cite{arcameron22,neumann,rystsov}).
In this paper, we prove three main results, showing different outcomes for
the ``low rank'' and the ``high rank'' cases. A graph admitting a vertex-primitive automorphism 
group will be called a \emph{primitive graph}.

\begin{theorem}
There are primitive graphs admitting non-uniform endomorphisms, and hence
not every primitive group is almost synchronizing.
\end{theorem}

Section \ref{smallrnk} describes a number of different graph constructions used to demonstrate this result.
This theorem resolves the almost synchronizing conjecture in the negative, but shows that the
structure of endomorphisms of primitive graphs can be more complex than
previously suspected, prompting the difficult problem of finding a classification of the 
primitive almost synchronizing groups.

Our second main theorem considers the high rank case and extends the
spectrum of ranks for which transformations of that rank are known to be synchronized by
every primitive group.
\begin{theorem}
A primitive group of degree $n$ synchronizes every transformation of rank $n-3$ or $n-4$.
\end{theorem}	

Sections \ref{largernk} and \ref{largernk2} detail the somewhat intricate arguments required to deal with the considerable number of graphs that arise
in proving this result.

Our third main theorem considers groups of small \emph{permutation rank}. Groups with permutation rank $2$ are doubly transitive
and easily seen to be synchronizing, so the first non-trivial case is for groups of rank $3$. In this case, the group acts primitively on a \emph{strongly regular} graph and properties of such graphs can be used to considerably extend the $n-4$ bound.

\begin{theorem}
A primitive permutation group of degree $n$ and permutation rank~$3$
synchronizes any non-permutation with rank at least $n - (1+\sqrt{n-1}/12)$.
\end{theorem}

The original context of this research is automata theory, as we now outline. Our automata are always finite and deterministic. On reading a symbol, an automaton undergoes a change of state; so each symbol defines a \emph{transition}, a transformation on the set of states. The set of transformations realised by reading a word or sequence of symbols is the semigroup generated by the transitions of the automaton. Thus, from an algebraic viewpoint, an automaton is a subsemigroup of the full transformation semigroup on a finite set with a prescribed set of generators.

A deterministic finite-state automaton is said to be \emph{synchronizing} if there is a finite word such that, after reading this word, the automaton is in a fixed state, independently of its state before reading the word. In other terms, the word evaluates to a transformation of rank~$1$, mapping the set of states to a single state. Such a word is called a \emph{reset word} or \emph{synchronizing word}.

One of the oldest and most famous problems in automata theory, the well-known  \v Cern\'y conjecture, states that if an automaton with $n$ states has a synchronizing word, then there exists one of length $(n-1)^2$. (For many references on the growing bibliography on this problem please see the two websites \cite{JEP,Tr} and also Volkov's talk \cite{Vo}; so far the best bound for the length of a reset word is cubic \cite{cubic}.)
Solving this conjecture is equivalent to proving that given a set $S$ of transformations on a finite set of size $n$ then, if
the transformation semigroup $\langle S \rangle$ contains a constant transformation, then it contains one that can be expressed
as a word of length at most $(n-1)^2$ in the generators of $S$. This conjecture has been established
for \emph{aperiodic automata}, that is, when $\langle S\rangle$ is a semigroup with no non-trivial subgroups \cite{Tr07}. So it remains to prove the conjecture for semigroups that do contain non-trivial subgroups, and the case when the semigroup contains a permutation group is a particular instance of this general problem. Indeed, the known examples witnessing the optimality of the \v Cern\'y bound contain a permutation among the given set of generators.


We note that if a transformation
semigroup $S$ contains a transitive group $G$ but not a constant function, then
the image $I$ of a transformation $f$ of minimum rank in $S$ is a
\emph{$G$-section} for the kernel of $f$, in the sense that $Ig$ is a \emph{section}
or \emph{transversal} for $\ker(f)$, a set meeting every kernel class
in a single element. In addition, the transformation $f$ has uniform kernel 
(see Neumann~\cite{neumann}).

Although they did not use this terminology, results due to Rystsov \cite{rystsov} and Neumann \cite{neumann} 
cover some cases of the almost sychronizing conjecture. In particular, Rystsov \cite{rystsov} showed that a
transitive permutation group of degree $n$ is primitive if and only if
it synchronizes every transformation of rank $n-1$, while Neumann \cite{neumann} showed
that a primitive permutation group synchronizes every transformation of rank $2$.

In earlier work, Ara\'{u}jo and Cameron \cite{arcameron22} resolved some
additional cases of the conjecture:
\begin{theorem}
A primitive permutation group $G$ of degree $n$ synchronizes maps of 
kernel type $(k,1,\ldots,1)$ (for $k\ge 2$) and maps of rank $n-2$, 
as well as non-uniform maps of rank $4$ or $3$.
\end{theorem}

That paper, like the present one, uses a graph-theoretic
approach due to the third author \cite{synch_ltcc}. In what follows, 
graphs are always simple and undirected. The \emph{clique number} of a graph
is the largest number of vertices in a complete subgraph, while the 
\emph{chromatic number} is the smallest number of colours required for a proper
colouring.

\begin{theorem}
A transformation semigroup does not
contain a constant transformation if and only if it is contained in the
endomorphism monoid of a non-null graph. Moreover, we may
assume that this graph has clique number equal to chromatic number.
\end{theorem}

One direction of the theorem is clear, since a non-null graph has no rank~$1$
endomorphisms. For the other direction, define a graph by
joining two vertices if no element of the semigroup maps them to the same
place, and show that this graph has the required property. We will elaborate
further in the next section.

\medskip

This paper is in three main parts. In the first part,  we show that the almost synchronizing conjecture fails for maps of small rank,
so the results in \cite{arcameron22} are best possible. We construct four
examples of primitive groups (with degrees $45$, $153$, $495$ and $495$)
which fail to synchronize non-uniform maps of rank~$5$.  In addition, we find infinitely
many examples for rank~$6$, along with yet another sporadic example of rank~$7$ of degree $880$. 
Also, we provide a construction of primitive graphs whose automorphism groups have approximately 
$\sqrt{n}$ \emph{non-synchronizing ranks}, refuting a conjecture of the third author's 
on the number of non-synchronising ranks of a primitive group.

In the second part, we press forward with maps of
large rank, showing that a primitive group synchronizes all maps with kernel
type $(p,2,1,\ldots,1)$ or kernel type $(p,3,1,\dots,1)$ for $p\ge3$, as well as
all maps of rank $n-3$ and $n-4$. We also show that a primitive group synchronizes every map in which one
non-singleton kernel class is sufficiently large compared to the other non-singleton kernel classes.

In the third part, we consider the special situation where the primitive group $G$ has
permutation rank $3$, in which case any graph with automorphism group containing $G$
is either trivial or strongly regular. In the latter case, we prove a general result about endomorphisms of
strongly regular graphs, and deduce that $G$ synchronizes every non-permutation
transformation of rank at least $n-(1+\sqrt{n-1}/12)$.

The paper ends with a number of natural but challenging problems related to
synchronization in primitive groups and in related combinatorial settings.

Given the enormous progress made in the last three or four decades, permutation groups now has the tools to answer questions coming from the real world through transformation semigroups; these questions translate into beautiful statements in the language of permutation groups and combinatorial structures, as shown in many recent investigations (as a small sample, please see \cite{anarca,abc,arcameron22,arcameron,acmn,ams,ArnoldSteinberg,randomsynch,gr,neumann,sv15}).

\section{Transformation semigroups and graphs}\label{trans}

The critical idea underlying our study is a graph
associated to a transformation semigroup in the following way.
If $S$ is a transformation semigroup on
$\Omega$, then form a graph, denoted $\Gr(S)$, with
vertex set $\Omega$ where $v$ and $w$ are adjacent if and only if there
is no element $f$ of $S$ which maps $v$ and $w$ to the same point.
Now the following result is almost immediate (\emph{cf.} \cite{CK,synch_ltcc}).

\begin{theorem}(See \cite{arcameron22,CK})\label{Sgraph}
Let $S$ be a transformation semigroup on $\Omega$ and let $\Gr(S)$ be defined as above.
Then
\begin{enumerate}
\item[(a)] $S$ contains a map of rank $1$ if and only if $\Gr(S)$ is null (i.e., edgeless).
\item[(b)] $S\le\End(\Gr(S))$, and $\Gr(\End(\Gr(S)))=\Gr(S)$.
\item[(c)] The clique number and chromatic number of $\Gr(S)$ are both equal to
the minimum rank of an element of $S$.
\end{enumerate}
\end{theorem}

In particular, if $S=\langle G,f\rangle$ for some group $G$, then
$G\le\Aut(\Gr(S))$. So, for example, if $G$ is primitive and does not
synchronize $f$, then $\Gr(S)$ is non-null and has a primitive automorphism group, and so is
connected.

In this situation, assume that $f$ is an element of minimum rank in $S$;
then the kernel of $f$ is a partition $\rho$ of $\Omega$, and its image $A$
is a \emph{$G$-section} for $\rho$ (that is, $Ag$ is a section for $\rho$,
for all $g\in G$). Neumann~\cite{neumann}, analysing this situation, defined
a graph $\Delta$ on $\Omega$ whose edges are the images under $G$ of the
pairs of vertices in the same $\rho$-class. Clearly $\Delta$ is a subgraph
of the complement of $\Gr(S)$, since edges in $\Delta$ can be collapsed by
elements of $S$. Sometimes, but not always, $\Delta$ is the complement of
$\Gr(S)$.

We now  introduce a refinement of the previous graph $\Gr(S)$, which will allow us to
obtain the results of the remaining cases more easily. The new graph is
denoted by $\Gr'(S)$.
The same construction was used in a different context in \cite{randomsynch},
where it was called the \emph{derived graph} of $\Gr(S)$.

Suppose that $\gr(S)$ has clique number and chromatic number $r$ (where $r$
is the minimum rank of an element of $S$). We define $\gr'(S)$ to be the graph
with the same vertex set as $\gr(S)$, and whose edges are all those edges of $\gr(S)$
which are contained in $r$-cliques of $\gr(S)$.

\begin{theorem}\label{th:new graph}
Let $S$ be a transformation semigroup on $\Omega$ and let $\Gr(S)$ and $\gr'(S)$ be defined as above.
Then
\begin{enumerate}
\item $S$ contains a map of rank $1$ if and only if $\gr'(S)$ is null.
\item $S\le\End(\gr(S))\le\End(\gr'(S))$.
\item The clique number and chromatic number of $\gr'(S)$ are both equal to
the minimum rank of an element of $S$.
\item Every edge of $\gr'(S)$ is contained in a maximum clique.
\item If $S=\langle G,f\rangle$, where $G$ is a primitive permutation group
and $f$ a map which is a non-permutation not synchronized by $G$,
then $\gr'(S)$ is neither complete nor null.
\end{enumerate}
\end{theorem}

\begin{pf}
Elements of $\End(\gr(S))$ preserve $\gr(S)$ and map maximum cliques to
maximum cliques, so $\End(\gr(S))\le\End(\gr'(S))$. The existence of an
$r$-clique and an $r$-colouring of $\gr'(S)$ are clear, and so (c) holds; then
(a) follows. Part (d) is clear from the definition. For (e), the hypotheses
guarantee that the minimum rank of an element of $S$ is neither $1$ nor $n$.
\qed
\end{pf}

Note that strict inequality can hold in (b). If $\Gamma$ is the disjoint
union of complete graphs of different sizes, then $\gr'(\End(\Gamma))$
consists only of the larger complete graph, and has more endomorphisms than
$\Gamma$ does.

The next lemma is proved in \cite{arcameron22}, but since the techniques it introduces are important in subsequent arguments we provide its proof here.

\begin{lemma}\label{neigh}
Let $X$ be a nontrivial graph and let $G\le \Aut(X)$ be primitive. Then no two vertices of $X$ can have the same neighbourhood.
\end{lemma}
\begin{pf}
For $a\in X$ denote its neighbourhood by $N(a)$. Suppose that $a,b\in X$, with $a\neq b$, and $N(a)=N(b)$. We are going to use two different techniques to prove that this leads to a contradiction. The first uses the fact that the graph has at least one edge, while the second uses the fact that the graph is not complete.

First.
Define the following relation on the vertices of the graph: for all $x,y\in X$,
\[
x\equiv y \Leftrightarrow N(x)=N(y).
\]
This is an equivalence relation and we claim that $\equiv$ is neither the universal relation nor the identity. The latter follows from the fact that by assumption $a$ and $b$ are different and $N(a)=N(b)$. Regarding the former, there exist adjacent vertices $c$ and $d$ (because $X$ is non-null);
now $c \in N(d)$
but $c \notin N(c)$, so $c \not\equiv d$.
As $G$ is a group of automorphisms of $X$ it follows that $G$ preserves $\equiv$, a non-trivial equivalence relation, and hence $G$ is imprimitive, a contradiction.

Second. Assume as above that we have $a,b\in X$ such that $N(a)=N(b)$. Then the transposition $(a\ b)$ is an automorphism of the graph. A primitive group containing a transposition is the symmetric group. (This well-known result from the
nineteenth century can be found, for example, in \cite[p.241]{isaacs}.)
Hence $X$ is the complete graph, a contradiction.
\qed
\end{pf}

We conclude this section recalling another result from \cite{arcameron22} about \emph{primitive graphs}.

\begin{lemma}(\cite{arcameron22})\label{primgr}
Let $\Gamma$ be a non-null graph with primitive automorphism group $G$, and
having chromatic number $r$. Then $\Gamma$ does not contain a subgraph
isomorphic to the complete graph on $r+1$ vertices with one edge removed.
\end{lemma}

This lemma was important in \cite{arcameron22} and it will be here too (please see the observations after Lemma \ref{lm:estimate}).

\section{Maps of small rank}\label{smallrnk}

In this section, we discuss various counterexamples to the conjecture that
primitive
groups are almost synchronizing. From the discussion above, it
suffices to find a non-null graph $\Gamma$ with a primitive
automorphism group, and then
exhibit a non-uniform proper endomorphism $f$ of $\Gamma$.
Such an endomorphism is then a witness
that $G$ is not almost synchronizing for any
primitive group $G \leq \mathrm{Aut}(\Gamma)$.

In the first subsection, we present a number of sporadic examples
of vertex-primitive graphs, each with non-uniform proper endomorphisms of
rank $5$ or $7$.
The smallest of these, on $45$ vertices,  can be shown (by computer) to be the
unique smallest counterexample to the almost-synchronizing conjecture; the
details are given in Section~\ref{compute}. None of the graphs described here
are Cayley graphs.

In the second subsection, we present infinite families of primitive
graphs with non-uniform proper endomorphisms of rank $6$ and above.

\subsection{Rank 5 (and 7)}

In this subsection, we give the first counterexample to the conjecture that
primitive groups are almost synchronizing. In particular, we
construct a primitive group of degree $45$ that fails to
synchronize a non-uniform map of rank~$5$ with kernel type $(5,5,10,10,15)$.
We give two proofs of this, which
extend in different ways.

Our primitive group is $\mathrm{P}\Gamma\mathrm{L}(2,9)$ (also known as $A_6:2^2$), acting on $45$ points (this is \texttt{PrimitiveGroup(45,3)} in both \textsf{GAP} and {\sc Magma}).
This group has a suborbit of length $4$, and the orbital graph $\Gamma$ has
the property that any edge is contained in a unique triangle: the
closed neighbourhood of a vertex is a ``butterfly'' consisting of two
triangles with a common vertex (see Figure~\ref{butterfly}). Indeed, this graph
is the line graph of the celebrated \emph{Tutte--Coxeter graph} on
$30$ vertices, which in turn is the incidence graph of the generalized quadrangle $W(2)$ of
order~$2$. The graph was first found by Tutte~\cite{tutte} with a geometric interpretation
by Coxeter~\cite{coxeter,tutte2}.

Let $D$ be a dihedral subgroup of order $10$ of the automorphism group of
the graph. It is clear that elements of
order $5$ in $D$ fix no vertices of the graph, and a little thought shows
that their cycles are independent sets in the graph.

The full automorphism group of the graph is the automorphism group of $S_6$
(that is, the group extended by its outer automorphism), and there are two
conjugacy classes of dihedral groups of order $10$. It is important to take
the right one here: we want the $D_{10}$ which is \emph{not} contained in
$S_6$.

For this group $D$, we find that each orbit of $D$ is an independent set in
$\Gamma$; so there is a homomorphism of $\Gamma$ in which each orbit is
collapsed to a single vertex. A small calculation shows that the image of
this homomorphism is the graph shown below:

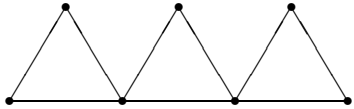
\begin{center}
\setlength{\unitlength}{1mm}
\begin{picture}(45,12.5)
\multiput(0,0)(15,0){4}{\circle*{1}}
\multiput(7.5,12.5)(15,0){3}{\circle*{1}}
\put(0,0){\line(1,0){45}}
\multiput(0,0)(15,0){3}{\line(3,5){7.5}}
\multiput(15,0)(15,0){3}{\line(-3,5){7.5}}
\end{picture}
\end{center}

Now this graph can be found as a subgraph of $\Gamma$, as the union of two
butterflies sharing a triangle; therefore the homomorphism can be realised
as an endomorphism of $\Gamma$ of rank~$7$, with kernel classes of sizes
$(10,10,5,5,5,5,5)$. The endomorphisms of ranks $5$ and $3$ can now be
 found by folding in one or
both ``wings'' in the above figure.

Our second approach uses the fact that the chromatic number and clique
number of this graph are each equal to $3$;  thus, each triangle has one
vertex in each of the three colour classes, each of size $15$. So there is
a uniform map of rank $3$ not synchronized by $G$.

We used \textsf{GAP} to construct the graph (the vertex numbering
is determined by the group), and its package \texttt{GRAPE} to find all the
independent sets of size $15$ in $\Gamma$ up to the action of $G$.
One of the two resulting sets is
\[A=\{ 1, 2, 3, 5, 10, 15, 16, 17, 25, 26, 27, 30, 42, 44, 45 \}.\]
The induced subgraph on the complement of this set has
two connected components, a $10$-cycle and a $20$-cycle. If we let $B$
and $C$ be the bipartite blocks in the $10$-cycle and $D$ and $E$ those
in the $20$-cycle, we see that $A,B,C,D,E$ are all independent sets, and
the edges between them are shown in Figure~\ref{butterfly}. Thus there is a proper endomorphism
mapping the graph to the closed neighbourhood of a vertex, with kernel classes $A,B,C,D,E$.

\begin{figure}
\begin{center}
\begin{tikzpicture}[scale=1.25]
\tikzstyle{vertex}=[circle,draw=black,fill=white,inner sep = 0.55mm, outer sep = 0.1mm]
\node [vertex] (v0) at (0,0) [label=above:{\footnotesize $A$}] {};
\node [vertex] (v1) at (0.87,0.5) [label=right:{\footnotesize $D$}]  {};
\node [vertex] (v2) at (0.87,-0.5)  [label=right:{\footnotesize $E$}]{};
\node [vertex] (v3) at (-0.87,-0.5)  [label=left:{\footnotesize $C$}] {};
\node [vertex] (v4) at (-0.87,0.5)  [label=left:{\footnotesize $B$}]{};
\draw (v0)--(v1)--(v2)--(v0)--(v3)--(v4)--(v0);
\end{tikzpicture}
\end{center}
\caption{The butterfly}\label{butterfly}
\end{figure}
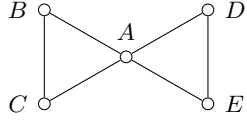

Using software developed at St~Andrews (see Section~\ref{compute} for details)
we were able to calculate all the proper endomorphisms of
this graph: there are $103680$ of these, with ranks $3$, $5$ and $7$;
the numbers of endomorphisms of each of these ranks are $25920$, $51840$ and
$25920$ respectively. Then, using GAP, we were able to determine
that the endomorphism monoid of this graph is given by $\End(X)=\langle G,t\rangle$,
where $G$ is $\mathrm{P\Gamma L}(2,9)$ and $t$ is the transformation

{\small
\begin{align*}
 t=&\text{Transformation}( [ 1, 1, 1, 14, 9, 14, 28, 41, 41, 1, 43, 28, 28, 41, 9, 1, 1, 25, 25, 28, 28, \\
 & 25, 41, 28, 1, 1, 9, 43, 14, 9, 43, 28, 28, 25, 41, 43, 14, 28, 43, 25, 14, 1, 28, 1, 9 ] ).
\end{align*}}

The endomorphisms of each possible rank form a single D-class.  The
structure for the H-classes is $S_3$, $D_8$ and $D_8$ for the three classes
respectively, where $D_8$ is the dihedral group on $4$ points. (Note that
these groups are the automorphism groups of the induced subgraphs on the
image of the maps.)

\medskip

A very similar example occurs in the line graph of the \emph{Biggs--Smith
graph} \cite{biggs, bs}, a graph on $153$ vertices whose automorphism group
is isomorphic to $\mathrm{PSL}(2,17)$ (\texttt{PrimitiveGroup(153,1)} in both
 \textsf{GAP} and {\sc Magma}). This graph has an endomorphism of rank $5$ and
kernel type $(6,6,45,45,51)$ constructed in a virtually identical way.

However, this particular construction gives no additional examples. A
vertex-primitive $4$-regular graph whose neighbourhood is
a butterfly is necessarily the linegraph of an {\em edge-}primitive
cubic graph. These were classified by Weiss \cite{weiss}, who
determined that the complete list is $K_{3,3}$, the Heawood
graph, the Tutte-Coxeter graph and the Biggs-Smith graph.
From either a direct analysis, or simply referring to
the small-case computations described in Section~\ref{compute},
it follows that the first two of these do not yield examples.

\medskip

However, we have found three additional examples with the help of the computer.
Surprisingly all three of them are associated with the group
$\mathrm{Aut}(M_{12})=M_{12}:2$. This group has two inequivalent primitive actions of
degree $495$. Each of them is the automorphism group of a graph of valency
$6$ in which the closed neighbourhood of a vertex consists of three triangles
with a common vertex, and in each case, the graph has chromatic number $3$.
In each case, there is a subgroup of the automorphism group with orbits of sizes
$55$, $55$, $110$, $110$, $165$; each orbit is an independent set and the connections
between the orbits give a homomorphism onto the butterfly.

The third example is associated with a {\em different} primitive action of $M_{12}:2$,
this time of degree $880$. In this action, $M_{12}:2$ is the full automorphism group of a $6$-regular graph
where each open neighbourhood is the disjoint union of two triangles. The group
has a subgroup of order $55$, which has $16$ equal-sized orbits each
inducing an independent set. These $16$ orbits can each be mapped to a
single vertex in such a way that the entire graph is mapped onto the closed neighbourhood of a vertex,
yielding an endomorphism of rank $7$, with kernel type $(220, 165, 165, 165, 55, 55, 55)$.
As the closed neighbourhood of a vertex consists of two $4$-cliques overlapping
in a vertex, we may perhaps view this just as a butterfly with bigger wings?

\subsection{Rank 6 and above}

While the constructions of the previous subsection seem to be sporadic
examples, we can also find several infinite families of vertex-primitive
graphs with proper non-uniform endomorphisms.

\begin{figure}
\begin{center}
\begin{tikzpicture}[scale=0.85]
\tikzstyle{vertex}=[circle,draw=black,fill=white,inner sep = 0.55mm, outer sep = 0.1mm]
\draw(0,0) grid (3,3);
\foreach \x in {0,1,2,3} {
\foreach \y in {0,1,2,3} {
\node [vertex] (v\x\y) at (\x, \y) {};
}}
\foreach \x in {0,1,2,3} {
\draw [bend left=25] (v\x0) to (v\x2);
\draw [bend left=25] (v\x1) to (v\x3);
\draw [bend right=15] (v\x0) to (v\x3);

\draw [bend left=25] (v0\x) to (v2\x);
\draw [bend left=25] (v1\x) to (v3\x);
\draw [bend right=15] (v0\x) to (v3\x);
}

\end{tikzpicture}
\end{center}
\caption{The graph $K_4 \cart K_4$}
\label{k4cartk4}
\end{figure}
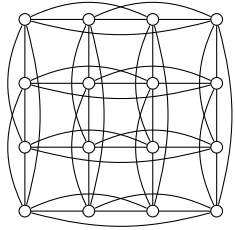

Recall that the {\em Cartesian product} $X \cart Y$ of two graphs $X$ and $Y$ is
the graph with vertex set $V(X \cart Y) = V(X) \times V(Y)$ and where
vertices $(x_1,y_1)$ and $(x_2, y_2)$ are adjacent if and only if they have equal
entries in one coordinate position and adjacent entries (in $X$ or $Y$ accordingly) in the other.
Figure~\ref{k4cartk4} shows the graph $K_4 \cart K_4$ both to illustrate the Cartesian product
and because it plays a role later in this section.

If $X$ is a vertex-primitive graph then the Cartesian product $X \cart X$ is also vertex-primitive, with automorphism group ${\rm Aut}(X) \wr \text{Sym(2)}$. In addition, if the chromatic and clique number of $X$ are both equal to $k$, then $V(X)$
can be partitioned into $k$ colour classes of equal size --- say $V_1$, $V_2$, $\ldots$, $V_k$, and there is a surjective homomorphism $X \cart X \rightarrow K_k \cart K_k$ with  kernel classes $\{V_{i} \times V_{j} \mid 1 \leq i,j \leq k\}$.
Therefore if there is a homomorphism $f: K_k \cart K_k \rightarrow X$, then by composing homomorphisms
\[
X \cart X \longrightarrow K_k \cart K_k \xlongrightarrow{f} X \longrightarrow X \cart X,
\]
there is an endomorphism of $X \cart X$. Moreover, if the homomorphism $f$ is non-uniform, then
the endomorphism is also non-uniform.

\begin{figure}
\begin{center}
\begin{tikzpicture}[scale=0.65]
\tikzstyle{vertex}=[circle,draw=white,fill=white,inner sep = 0.25mm, outer sep = 0.1mm]

\foreach \x in {0,1,2,3} {
  \draw  (\x,-0.5)--(\x,3.5);
}
\foreach \x in {0,1,2,3} {
  \draw  (-0.5,\x)--(3.5,\x);
}

\node [vertex] (v00) at (0,0) {\scriptsize {\tt 00}};
\node [vertex] (v01) at (0,1) {\scriptsize {\tt 11}};
\node [vertex] (v02) at (0,2) {\scriptsize {\tt 22}};
\node [vertex] (v03) at (0,3) {\scriptsize {\tt 33}};
\node [vertex] (v10) at (1,0) {\scriptsize {\tt 11}};
\node [vertex] (v11) at (1,1) {\scriptsize {\tt 00}};
\node [vertex] (v12) at (1,2) {\scriptsize {\tt 33}};
\node [vertex] (v13) at (1,3) {\scriptsize {\tt 22}};
\node [vertex] (v20) at (2,0) {\scriptsize {\tt 22}};
\node [vertex] (v21) at (2,1) {\scriptsize {\tt 33}};
\node [vertex] (v22) at (2,2) {\scriptsize {\tt 01}};
\node [vertex] (v23) at (2,3) {\scriptsize {\tt 10}};
\node [vertex] (v30) at (3,0) {\scriptsize {\tt 33}};
\node [vertex] (v31) at (3,1) {\scriptsize {\tt 22}};
\node [vertex] (v32) at (3,2) {\scriptsize {\tt 10}};
\node [vertex] (v33) at (3,3) {\scriptsize {\tt 01}};

\pgftransformxshift{6cm}

\foreach \x in {0,1,2,3} {
  \draw  (\x,-0.5)--(\x,3.5);
}
\foreach \x in {0,1,2,3} {
  \draw  (-0.5,\x)--(3.5,\x);
}

\node [vertex] (v00) at (0,0) {\scriptsize {\tt 00}};
\node [vertex] (v01) at (0,1) {\scriptsize {\tt 11}};
\node [vertex] (v02) at (0,2) {\scriptsize {\tt 22}};
\node [vertex] (v03) at (0,3) {\scriptsize {\tt 33}};
\node [vertex] (v10) at (1,0) {\scriptsize {\tt 11}};
\node [vertex] (v11) at (1,1) {\scriptsize {\tt 02}};
\node [vertex] (v12) at (1,2) {\scriptsize {\tt 33}};
\node [vertex] (v13) at (1,3) {\scriptsize {\tt 20}};
\node [vertex] (v20) at (2,0) {\scriptsize {\tt 22}};
\node [vertex] (v21) at (2,1) {\scriptsize {\tt 33}};
\node [vertex] (v22) at (2,2) {\scriptsize {\tt 10}};
\node [vertex] (v23) at (2,3) {\scriptsize {\tt 01}};
\node [vertex] (v30) at (3,0) {\scriptsize {\tt 33}};
\node [vertex] (v31) at (3,1) {\scriptsize {\tt 20}};
\node [vertex] (v32) at (3,2) {\scriptsize {\tt 01}};
\node [vertex] (v33) at (3,3) {\scriptsize {\tt 12}};

\pgftransformxshift{6cm}

\foreach \x in {0,1,2,3} {
  \draw  (\x,-0.5)--(\x,3.5);
}
\foreach \x in {0,1,2,3} {
  \draw  (-0.5,\x)--(3.5,\x);
}

\node [vertex] (v00) at (0,0) {\scriptsize {\tt 00}};
\node [vertex] (v01) at (0,1) {\scriptsize {\tt 11}};
\node [vertex] (v02) at (0,2) {\scriptsize {\tt 22}};
\node [vertex] (v03) at (0,3) {\scriptsize {\tt 33}};
\node [vertex] (v10) at (1,0) {\scriptsize {\tt 11}};
\node [vertex] (v11) at (1,1) {\scriptsize {\tt 02}};
\node [vertex] (v12) at (1,2) {\scriptsize {\tt 33}};
\node [vertex] (v13) at (1,3) {\scriptsize {\tt 20}};
\node [vertex] (v20) at (2,0) {\scriptsize {\tt 23}};
\node [vertex] (v21) at (2,1) {\scriptsize {\tt 30}};
\node [vertex] (v22) at (2,2) {\scriptsize {\tt 01}};
\node [vertex] (v23) at (2,3) {\scriptsize {\tt 12}};
\node [vertex] (v30) at (3,0) {\scriptsize {\tt 32}};
\node [vertex] (v31) at (3,1) {\scriptsize {\tt 23}};
\node [vertex] (v32) at (3,2) {\scriptsize {\tt 10}};
\node [vertex] (v33) at (3,3) {\scriptsize {\tt 01}};

\end{tikzpicture}
\end{center}
\caption{Non-uniform homomorphisms from $K_4 \cart K_4$ to its complement}
\label{fig:homs}
\end{figure}
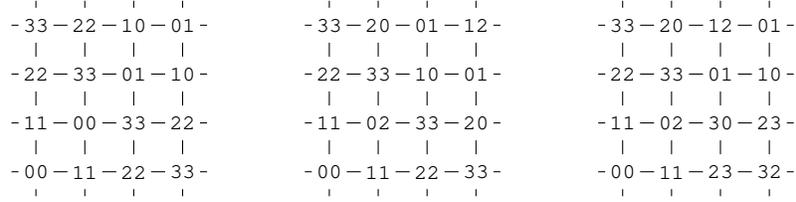

Although at first sight, there appears to be little to be gained from this observation, in practice it is much easier (both computationally and theoretically) to find homomorphisms between the two relatively small graphs $K_k \cart K_k$ and $X$, than working directly with the larger graph $X \cart X$. Although this finds only a restricted subset of endomorphisms, it turns out to be sufficient to find large numbers of non-uniform examples.

We start by considering the case where $X = \overline{K_k \cart K_k}$, where we
assume that the vertices of both graphs are labelled with pairs $(i,j)$, where $0 \leq i,j < k$ and that in $K_k \cart K_k$ two distinct vertices are adjacent if and only if they agree in one coordinate position, while in its complement, they are adjacent if and only if they {\em disagree} in both coordinate positions.

The graph $X$ has chromatic number and clique number equal to $k$. (A diagonal
set $\{(i,i):0\le i< k\}$ is a $k$-clique, while using the first coordinate as
colour gives a $k$-colouring.) The homomorphisms we seek are those
from $K_k \cart K_k$ to its own complement.

In particular, Figure~\ref{fig:homs} exhibits three non-uniform homomorphisms of ranks $6$, $9$ and $12$ from $K_4 \cart K_4$ to its complement, where the diagrams show the image of each vertex, but using $xy$ to represent $(x,y)$.
Verifying that this function is a homomorphism merely requires checking for each row that the four pairs assigned to it have pairwise distinct first co-ordinates and pairwise distinct second co-ordinates, and similarly for each column. This happens if and only if the pairs are obtained from the super-position of two Latin squares of order $4$, one determining the first co-ordinate and the other the second co-ordinate. The rank of the homomorphism is then just the total number of distinct pairs that occur --- this number ranges from a minimum of $k$ (when the two Latin squares are identical) to a maximum of $k^2$ (when the two Latin squares are orthogonal). In the example of Figure~\ref{fig:homs} the kernel types of the homomorphisms  are $\{2^4, 4^2\}$,  $\{1^4, 2^4, 4\}$ and $\{1^8, 2^4\}$, which correspond to endomorphisms of $X \cart X$ of the same rank, but with kernel classes each $16$ times larger.

This argument clearly generalises to all $k \geq 4$ (non-uniform homomorphisms do not arise when $k < 4$) and so any two Latin squares of order $k$ ({\em not} necessarily orthogonal) will determine an endomorphism of $\overline{K_k \cart K_k}$. Two Latin squares are said to be {\em $r$-orthogonal} if $r$ distinct pairs arise when they are superimposed. Thus we find an endomorphism of rank $r$ from any pair of $r$-orthogonal Latin squares.  The following result, due to Colbourn \& Zhu \cite{colbournzhu} and Zhu \& Zhang \cite{zhuzhang} shows exactly which possible ranks arise in this fashion.

\begin{theorem}
There are two $r$-orthogonal Latin squares of order $k$ if and only if $r \in \{k, k^2\}$ or $k+2 \leq r \leq k^2-2$, with the
following exceptions:

\begin{enumerate}
\item $k=2$ and $r=4$;
\item $k=3$ and $r \in \{5,6,7\}$;
\item $k=4$ and $r \in \{7,10,11,13,14\}$;
\item $k=5$ and $r \in \{8,9,20,22,23\}$;
\item $k=6$ and $r \in \{33,36\}$. \qed
\end{enumerate}
\end{theorem}

In particular, for any $k \geq 4$, there is an endomorphism of rank $k+2$ with image two $k$-cliques overlapping in
a $(k-2)$-clique. As $k$ increases, we get a sequence of butterflies  with increasingly fat bodies, but fixed-size wings.

This construction also sheds some light on the possible non-synchronizing ranks
for a group. For a group $G$ of degree $n$, a \emph{non-synchronizing rank}
is a value $r$ satisfying $2 \leq r \leq n-1$ such that $G$ fails to
synchronize some transformation of rank $r$.

A transitive imprimitive group of degree $n$, having $m$ blocks of
imprimitivity each of size $k$, preserves both a disjoint union of $m$
complete graphs of size $k$ (which has endomorphisms of ranks all multiples
of $k$) and the complete $m$-partite graph with parts of size $k$ (which has
endomorphisms of all ranks between $k$ and $n$ inclusive). From this, a short
argument shows that such a group has at least $(3/4 - o(1))n$ non-synchronizing
ranks. It was suspected that a primitive group has many fewer non-synchronizing
ranks, perhaps as few as $O(\log n)$. However, as this construction provides
approximately $k^2$ non-synchronizing ranks for a group of degree $k^4$, this
cannot be the case.

It is natural to wonder whether this construction can be used to find non-uniform endomorphisms for graphs {\em other than} $X = \overline{K_k \times K_k}$.
Unsurprisingly, the answer to this question is yes, with the line graph of the complete graph $L(K_n)$ (also known as the {\em triangular graph}) being a suitable candidate for $X$ whenever $n$ is even.  For example $L(K_6)$ is a $15$-vertex graph with chromatic number and clique number equal to $5$. The vertices of $L(K_6)$ can be identified with the endpoints of the corresponding edge in $K_6$, and thus each vertex of $L(K_6)$ is represented by a $2$-set of the form $\{x,y\}$ which we will abbreviate to $xy$.

\begin{figure}
\begin{center}
\begin{tikzpicture}[scale=0.75]
\tikzstyle{vertex}=[circle,draw=white,fill=white,inner sep = 0.25mm, outer sep = 0.1mm]

\foreach \x in {0,1,2,3,4} {
  \draw  (\x,-0.5)--(\x,4.5);
}
\foreach \x in {0,1,2,3,4} {
  \draw  (-0.5,\x)--(4.5,\x);
}

\node [vertex] (v00) at (0,0) {\footnotesize {\tt 01}};
\node [vertex] (v11) at (1,1) {\footnotesize {\tt 02}};
\node [vertex] (v22) at (2,2) {\footnotesize {\tt 03}};
\node [vertex] (v33) at (3,3) {\footnotesize {\tt 04}};
\node [vertex] (v44) at (4,4) {\footnotesize {\tt 05}};

\node [vertex] (v10) at (1,0) {\footnotesize {\tt 12}};
\node [vertex] (v20) at (2,0) {\footnotesize {\tt 13}};
\node [vertex] (v30) at (3,0) {\footnotesize {\tt 14}};
\node [vertex] (v40) at (4,0) {\footnotesize {\tt 15}};

\node [vertex] (v01) at (0,1) {\footnotesize {\tt 12}};
\node [vertex] (v21) at (2,1) {\footnotesize {\tt 23}};
\node [vertex] (v31) at (3,1) {\footnotesize {\tt 24}};
\node [vertex] (v41) at (4,1) {\footnotesize {\tt 25}};

\node [vertex] (v02) at (0,2) {\footnotesize {\tt 13}};
\node [vertex] (v12) at (1,2) {\footnotesize {\tt 23}};
\node [vertex] (v32) at (3,2) {\footnotesize {\tt 34}};
\node [vertex] (v42) at (4,2) {\footnotesize {\tt 35}};

\node [vertex] (v03) at (0,3) {\footnotesize {\tt 14}};
\node [vertex] (v13) at (1,3) {\footnotesize {\tt 24}};
\node [vertex] (v23) at (2,3) {\footnotesize {\tt 34}};
\node [vertex] (v43) at (4,3) {\footnotesize {\tt 45}};

\node [vertex] (v04) at (0,4) {\footnotesize {\tt 15}};
\node [vertex] (v14) at (1,4) {\footnotesize {\tt 25}};
\node [vertex] (v24) at (2,4) {\footnotesize {\tt 35}};
\node [vertex] (v34) at (3,4) {\footnotesize {\tt 45}};

\end{tikzpicture}
\end{center}
\caption{A homomorphism from $K_5 \cart K_5$ to $L(K_6)$}
\label{fig:lk6}
\end{figure}
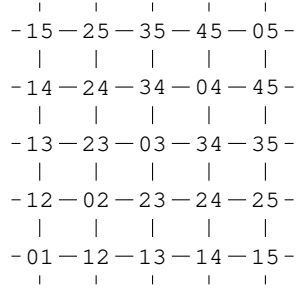

Figure~\ref{fig:lk6} depicts a surjective homomorphism from $K_5 \cart K_5$ to $L(K_6)$ by labelling each of the vertices of $K_5 \cart K_5$ with its image in $L(K_6)$  under the homomorphism. For each of the horizontal or vertical lines --- corresponding to the cliques of $K_5 \cart K_5$ --- it is easy to confirm that the images of the five vertices in the line share a common element and thus are mapped a clique of $L(K_6)$. This homomorphism has rank $15$ and kernel type $\{1^5, 2^{10}\}$ and hence yields a non-uniform endomorphism of $L(K_6) \cart L(K_6)$ with kernel classes of $25$ times the size. The pattern shown in Figure~\ref{fig:lk6} can be generalised to all triangular graphs, by defining a map  $f: K_{n-1} \cart K_{n-1} \longrightarrow L(K_n)$ by $f\left( (a,b) \right) = \{a+1,b+1\}$ if $a \not= b$ and $f \left( (a,a) \right)= \{0,a+1\}$. This homomorphism has kernel type $\{1^{n-1}, 2^{(n-1)(n-2)/2}\}$.

We finish this section with yet another construction that provides an infinite family of rank $6$ non-uniform non-synchronizable transformations.

Let $p$ be a prime greater than $5$, and let $V$ be the vector space spanned
by $e_0,\ldots,e_{p-1}$ (we think of the indices as elements of the integers mod $p$)
with the single relation that their sum is zero. Let $\Gamma$ be the Cayley graph for $V$
with connection set of size $2p$ consisting of the vectors $e_i$ and
$e_i+e_{i+1}$ with $i$ running over the integers mod $p$. It is clear that
the group $V:D_{2p}$ acts as automorphisms of this graph, and is primitive
provided that $2$ is a primitive root mod~$p$ (this is the condition for $V$
to be irreducible as a $C_p$-module).

Now let $X$ be the subspace spanned by $e_i+e_{i+2}$ for $i=0,1,\ldots,
p-5$. These vectors are linearly independent and so span a space of
codimension $3$, with $8$ cosets. Check that this subspace contains no
edge of the graph: no two of its vectors differ by a single basis vector
or a sum of two consecutive basis vectors. We can take coset representatives to
be $0,e_0,e_1,e_0+e_1,e_{p-2},e_{p-2}+e_0,e_{p-2}+e_1,e_{p-2}+e_0+e_1$.

Each coset contains no edges of the graph, and indeed the unions
$(X+e_1)\cup(X+e_{p-2})$ and $(X+e_0+e_1)\cup(X+e_0+e_{p-2})$ also contain
no edges. Collapsing these two unions and the other four cosets to a vertex,
we find by inspection that the graph is a ``butterfly'':

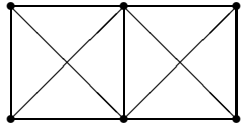
\begin{center}
\setlength{\unitlength}{1mm}
\begin{picture}(30,15)
\multiput(0,0)(15,0){3}{\circle*{1}}
\multiput(0,15)(15,0){3}{\circle*{1}}
\multiput(0,0)(0,15){2}{\line(1,0){30}}
\multiput(0,0)(15,0){3}{\line(0,1){15}}
\multiput(0,0)(15,0){2}{\line(1,1){15}}
\multiput(0,15)(15,0){2}{\line(1,-1){15}}
\end{picture}
\end{center}

The two vertices forming the butterfly's body are the ``double cosets''.

Now we can find a copy of the butterfly in the graph, using the vertices
$0$ and $e_0$ for the body, $e_1$ and $e_0+e_1$ for one wing, and
$e_{p-1}$ and $e_0+e_{p-1}$ for the other wing.

So there is an endomorphism of rank $6$, with two kernel classes of size
$2^{p-3}$ and four of size $2^{p-4}$.

\section{Maps of rank $n-3$}\label{largernk}

The goal of this section is to show that primitive groups synchronize maps of
rank $n-3$. Moreover, we will establish general properties of graphs with
primitive automorphism groups that will also be used in the next
sections; therefore these results will be stated and proved as generally as
possible.

\subsection{Groups with elements of small support}
\label{smallmindegree}

At several points during the argument (in this and in the next sections), we will establish that the (primitive)
automorphism group of a graph under consideration contains
a permutation which is a product of three or four disjoint transpositions,
and hence has support of size $m\in\{6,8\}$.
As the automorphism group of a non-trivial graph is not $2$-transitive (in particular, it does not contain $A_n$),
it follows from \cite{LS} that it has degree less than $(m/2+1)^2$. Thus
in the worst case, we need only examine primitive groups of degree less than
$25$, which can easily be handled computationally.
For convenience all of the computational results
are collated in Section~\ref{compute}.

\subsection{A bound on the intersection of neighbourhoods}\label{bounds}

The goal of this subsection is to prove the following result about $\gr(S)$, the graph introduced in Section \ref{trans}. (We note that Spiga and Verret \cite{sv15} have proved some results about neighbourhoods of vertex-primitive graphs that also imply this property.)

\begin{theorem}\label{k-2}
Let $G$ be a group acting primitively on a set $X$, and suppose that $f \in T(X)$ is not synchronized by $G$.
Let $S$ be the semigroup generated by $G$ and $f$, and let $k$ be the valency
of the graph $\Gam=\gr(S)$. Then for all distinct $x,y \in X$, their neighbourhoods $N(x), N(y)$ in $\Gam$ satisfy $|N(x) \cap N(y)|\le k-2$.
\end{theorem}

The previous theorem is an important ingredient in the proof of the main results in this paper. Its proof is based on
 methods from \cite{abc,arcameron22,CK} and is a consequence of the following sequence of lemmas.

\begin{lemma}\label{clique}
Let $\Gamma$ be a non-null graph with primitive automorphism group $G$, and
having chromatic number and clique number $r$. Let $x$ be a vertex of $\Gamma$, and $C$ an $r$-clique
containing $x$. Then for every vertex $y\not\in N(x) \cup\{ x\}$, we have that $(N(x)\setminus N(y)) \cap C \ne \emptyset.$
\end{lemma}
\begin{pf}
Assume instead that for some $y\not\in N(x)\cup\{x\}$ we have  $(N(x)\setminus N(y)) \cap C = \emptyset.$
Then every element of $C$, different from $x$, is a neighbour of $y$. Thus the set $C \cup \{y\}$ induces a subgraph that is isomorphic to the
complete graph with one edge removed. This contradicts Lemma \ref{primgr}.
\qed
\end{pf}

We next  state an observation on primitive groups and
quasiorders (reflexive and transitive relations). The proof
is an easy exercise.

\begin{lemma}\label{quasi}
Let $G$ be a permutation group on the finite set $X$. Then $G$ is primitive if and
only if the only $G$-invariant quasiorders are the identity and the universal
relation.
\end{lemma}
%
%
%

This immediately implies:
\begin{lemma}\label{quasiord}
Let $\Gam$ be a graph with primitive automorphism group $G$ and clique number $r$ on the vertex set $X$.
Assume that there are distinct elements
$a, b\in X$ satisfying the following property: every $r$-clique containing $a$ also contains $b$. Then $\Gam$ is complete.
\end{lemma}

\begin{pf}
The relation $\rightarrow$, defined by
$x\rightarrow y$ if every $r$-clique containing $x$ also contains $y$, is easily seen
to be a quasiorder, and so the result follows from Lemma~\ref{quasiord}.
%
\end{pf}

The following corollary follows from  Lemma \ref{quasiord}, Lemma \ref{clique}, and Lemma \ref{neigh}.
\begin{cor}\label{cork-2}
Let $\Gamma$ be a non-complete, non-null graph with primitive automorphism group $G$, 
clique number $r$ equal to its chromatic number, and valency $k$.
Then for any two distinct vertices $x$, $y$ we have  $|N(x) \cap N(y)| \le k-2$.
\end{cor}

\begin{pf}
By Lemma \ref{neigh}, no two vertices of $\Gam$ have the same neighbourhood. Hence it suffices to show that there are no distinct vertices
$x$, $y$ satisfying $|N(x) \cap N(y)| = k-1$.  For a contradiction, suppose that $x$ and $y$ have this property.

Assume first that  $x$ and $y$ are not adjacent, and let $z$ be the unique element in $N(x) \setminus N(y)$. Now let $C$ be a
clique of size $r$ containing $x$. As $y\not\in N(x) \cup \{x\}$, it follows from Lemma \ref{clique} that  $(N(x) \setminus N(y)) \cap C \not= \emptyset$, and as $N(x) \setminus N(y) = \{z\}$, it follows that $z \in C$. Therefore every $r$-clique containing $x$ also contains $z$, and
so by Lemma~\ref{quasiord}, $\Gam$ is complete, contradicting the hypotheses on $\Gam$.

If on the other hand $x$ and $y$ are adjacent, then
$N(x)\cup\{x\}=N(y)\cup\{y\}$, which defines a non-trivial $G$-invariant
equivalence relation on $X$, and so $\Gam$ is empty or complete, once again
contradicting the hypotheses on $\Gam$.

We conclude that $|N(x) \cap N(y)| \not= k-1$, and so $|N(x) \cap N(y)| \le k-2$.\qed
\end{pf}

Theorem \ref{k-2} is now an immediate consequence of the previous corollary and Theorem \ref{Sgraph}.
The results above apply in particular to the graph $\gr'(S)$, where
$S=\langle G,f\rangle$, since the automorphism group of this graph contains the primitive group $G$.

\begin{prop}\label{pr:derived}
Suppose that $G$ is primitive and does not synchronize $f$. Let
$S = \langle G,f \rangle$. Then $\gr'(S)$ has the following properties:
\begin{enumerate}
\item If $x\ne y$, then $|N(x)\cap N(y)|\le k-2$, where $k$ is the valency of
$ \gr'(S)$.
\item If $x$ and $y$ are distinct, there exists  a maximum clique in $\gr'(S)$
containing $x$ but not $y$.
\end{enumerate}
\end{prop}

\begin{pf}
The first claim follows  from Lemma \ref{quasiord} and Corollary \ref{cork-2}.
The second is clear since $\gr'(S)$ has the same maximum cliques as $\gr(S)$.
\qed
\end{pf}

\vspace{0.2cm}

\subsection{The main result about maps of rank $n-3$}\label{main section}

The aim of this subsection is to prove the following:

\begin{theorem}\label{main}	
Primitive groups synchronize maps of rank $n-3$.
\end{theorem}	

It is proved in \cite{arcameron22} that a primitive group $G$ synchronizes every map of kernel type $(4,1,\dots,1)$.
Therefore, to cover all maps of rank $n-3$, we have to consider the maps of kernel type $(3,2,1, \dots, 1)$
and $(2,2,2,1,\dots,1)$.

\subsubsection*{Kernel type $(p,2,1,\dots,1)$} 

It was shown in \cite{arcameron22} that every primitive group $G$ synchronizes every map $f$ of kernel type $(p,2,1,\dots,1)$ for $p=2$, and for idempotent maps in the case of $p=3$.
We will show that every primitive group $G$ synchronizes every map $f$ of kernel type $(p,2,1,\dots,1)$,
for $p\ge 3$.
 This result was recently proved
independently by Spiga and Verret \cite{sv15}.

\begin{theorem}\label{kerp211}
Let $X$ be a set with at least $6$ elements, $p \ge 3$,  $G$ a primitive group acting on $X$ and $f \in T(X)$ a map of kernel type $(p,2,1, \dots, 1)$, that is, $f$ has one kernel class of size $p$, one kernel class of size $2$, and an arbitrary number
of singleton kernel classes. Then $G$ synchronizes $f$.
\end{theorem}

\begin{pf}
Let $S= \langle G,f\rangle$, and $\Gam=\gr(S)$, let $k$ be the valency of $\Gam$.
 Assume that $G$ does not synchronize
$f$; then $\Gam$ is not null by Theorem \ref{Sgraph}, and it is not complete either as $f$ has
non-singleton kernel classes.

Let $A=\{a_1, a_2\}$ be the two-element kernel class of $f$ and let $B$ be its largest
 kernel class. Let $b_1,b_2$ be distinct elements in $B$, and $K =A \cup B$.

Now let  $N_B$ be the set of all vertices in $\overline{K}$, the complement of $K$, that are adjacent to at least one element of $B$.
As $f$ maps $N_B$ injectively into $N(b_1 f) \setminus \{a_1 f\}$, we get  $|N_B| \le k-1$.

By definition we have $N(b_1), N(b_2) \subseteq N_B \cup A$, and $$|N_B \cup A|= |N_B|+|A|\leq (k-1)+2=k+1.$$ As
$|N(b_1)|=|N(b_2)|= k$, it follows that $|N(b_1) \cap N(b_2)| \ge k-1$ by the pigeonhole principle.
This contradicts Theorem \ref{k-2}, and so $G$ synchronizes $f$.\qed
\end{pf}

\subsubsection*{Kernel type $(2,2,2,1,\dots,1)$} 

We are going to prove the following result:
\begin{theorem}\label{2221}
Let $G$ act primitively on $X$ and let $f \in T(X)$ have kernel type  $(2,2,2,1,\dots,1)$.
Then $G$ synchronizes $f$.
\end{theorem}

Let $S=\langle G,f\rangle$ and $\Gam =\gr(S)$. By Theorem \ref{Sgraph}, $S$ is a set of endomorphisms of $\Gam$.
Assume that $G$ does not synchronize $f$; then $\Gam$ is not null, once again by Theorem \ref{Sgraph}. Moreover $\Gam$  has clique number equal to its chromatic number.
Let $k$ be the valency of $\Gamma$.

Let $A=\{a_1, a_2\}$, $B=\{b_1, b_2\}$, $C=\{c_1,c_2\}$ be the  non-singleton kernel classes of $f$, and let $K=A \cup B \cup C$.

By \cite{neumann}, the smallest non-synchronizing group has degree $9$; therefore, every primitive group of degree at most $8$ synchronizes every singular transformation; hence we can assume that $n\ge 9$ and so
$f$ of kernel type $(2,2,2,1,\ldots)$ has  at least $3$ singletons classes so that $\overline{K} \ne \emptyset$.
As $\Gam$ has primitive automorphism group, it is connected and hence there is at least
one edge between $K$ and $\overline{K}$, say at some $a_i \in A$. We claim that there
is an edge between $A$ and $B \cup C$. For the sake of contradiction,  assume otherwise.  Then, as $f$ maps $\overline{K}$ injectively, both  $N(a_1)$ and $N(a_2)$ are mapped injectively to $N(a_1f)$ and as all of these sets have size $k$, we get that $N(a_1)=N(a_2)$, contradicting Lemma  \ref{neigh}.
So there is an edge between $A$ and, say $B$, and hence between $a_1 f $ and $b_1 f$.

Repeating the same argument for the remaining class $C$ we get that there must also be at least one edge between $C$ and one of $A$ or $B$. Up to a renaming of the
classes, we have two situations:

\subparagraph {Case 1:} there are no edges between $A$ and $C$.

\medskip
We exclude this case with an argument already used. For there are no edges between $A$ and $C$,
and so any neighbour of $a_1$ or $a_2$ must lie in $B \cup N_A$, where $N_A  $ is the set of
elements in $\overline{K}$ that is adjacent to at least one of $a_1,a_2$.

Now $|N_A| \le k-1$, as its elements are mapped injectively to $N(a_1 f)\setminus \{b_1 f\}$
by $f$, and so $|B \cup N_A| \le k+1$. By the pigeonhole principle, $|N(a_1) \cap N(a_2)| \le k-1$,
contradicting Theorem \ref{k-2}.

\subparagraph {Case 2:} there are edges between every pair from $A$, $B$, and $C$, and hence their images $a_1 f, b_1 f, c_1 f$ form a 3-cycle.

\medskip

Consider the induced subgraph on $X f$.
 We will obtain upper and lower bounds on the number of edges
in $X f$, using methods analogous to those used  in \cite{arcameron22}.

Let $e$ be the number of vertices  in $\Gam$, and let $l$ be the number of edges within $K$.

As $X f$ is obtained by deleting  three vertices of $X$, the induced graph on $X f$ contains
at most $e-3k +3$ edges (a loss of $k$ edges at each  vertex not in the image of $f$,
with at most
$3$ edges counted twice).

For the lower bound we count how many edges are at most
sent to a common image
 by $f$. Let $r,s,t$ be the number of edges between $A$ and $B$, $B$ and $C$, $C$ and $A$, respectively, hence $l=r+s+t$. Since the sets $A$, $B$ and $C$ each have two elements, it follows that  $r,s,t \le 4$.
These $l$ edges are collapsed onto $3$ edges, so we loose $l-3$ edges from within $K$.

For each $c \in \overline{K}$, such that $(c,a_1), (c,a_2)$ are edges, we map two edges into one (and hence lose one).
Now $r+t$ edges connect $A$ to $K \setminus A$. These edges are connecting to just two vertices, namely $a_1$ and $a_2$, so one of them connects
 to at least $\lceil (r+t)/2\rceil$ edges from within $K$.
Hence there are at most $k- \lceil (r+t)/2\rceil$ edges between one of the
 $a_i$ and $\overline{K}$, and so this is the maximal number of values $c \in \overline{K}$ for which $(c,a_1)$ and $(c,a_2)$ are edges. Symmetric arguments yield the following result.
\begin{lemma}\label{collapse2} The transformation
$f$ identifies at most $k- \lceil(r+t)/2\rceil$ of the edges between $\overline{K}$ and $A$, at most $k- \lceil (r+s)/2\rceil$ of the edges between $\overline{K}$ and $B$, and at most
$k- \lceil (s+t)/2\rceil$ of the edges between $\overline{K}$ and $C$.
\end{lemma}

Hence the number of edges in $X f$ is at least
$$e - \underbrace{(k-\lceil (r+t)/2\rceil)}_{\mbox{loss in $\overline{K}$--$A$}}-  \underbrace{(k-\lceil (r+s)/2\rceil)}_{\mbox{loss in $\overline{K}$--$B$}}
- \underbrace{(k-\lceil (s+t)/2 \rceil)}_{\mbox{loss in $\overline{K}$--$C$}} - \underbrace{(l-3)}_{\mbox{loss within $K$}} \ge$$
\begin{equation}\ge  e-3k + (r+t)/2+ (r+s)/2+ (s+t)/2 - (l-3) = e -3k +3.
\label{ineq}\end{equation}
This equals our upper bound. It follows that all estimates used in deriving our bounds must be tight. We have proved half of the following result.

\begin{lemma}\label{Case1}
Under the conditions of Case 2, and with notation as above, $f$ identifies exactly $k-  (r+t)/2$ pairs of  edges between $\overline{K}$ and $A$,  $k- (r+s)/2 $ pairs of  edges between $\overline{K}$ and $B$, and
$k- (s+t)/2 $ pairs of the edges between $\overline{K}$ and $C$. In addition,
$N(a_1) \cap \overline{K} = N(a_2) \cap\overline{K}$,  $N(b_1) \cap \overline{K} = N(b_2)\cap \overline{K}$,  $N(c_1) \cap \overline{K} = N(c_2) \cap \overline{K}$.
\end{lemma}
\begin{pf}
As our upper and lower bounds agree, the estimates from Lemma \ref{collapse2} must be tight.
Moreover,
the inequality in (\ref{ineq}) must be tight, as well,
 which implies that $(r+t)/2, (r+s)/2, (s+t)/2$ are equal to their ceilings and hence integers. This proves the first claim.

Thus  $2k-(r+t)$, the number of edges between $A$ and $\overline{K}$, is an even number. In addition, $k- (r+t)/2$, the number of edges
between $A$ and $\overline{K}$ identified by $f$, is then exactly half of the number of edges between $A$ and $\overline{K}$. However $f$ can only map at most two such edges onto one, as $A$ has only $2$ elements.
It follows that if $(c,a_1)$ is an edge with $c \in \overline{K}$, then $(c,a_2)$ is an edge as well, and vice versa. Hence $N(a_1) \cap \overline{K} = N(a_2) \cap\overline{K}$, and the remaining claims follow by symmetry.
\phantom{*}\hfill\qed
\end{pf}

Theorem \ref{k-2} implies that $N(b_1) \cup N(b_2)$ must contain at least four vertices that
are in exactly one of $N(b_1), N(b_2)$.
By Lemma \ref{Case1}, $N(b_1) \cap \overline{K} = N(b_2)\cap \overline{K}$. So the four vertices
that are in exactly one of $N(b_1), N(b_2)$  must be $a_1,a_2,c_1,c_2$,
with each of $b_1$ and $b_2$ connected
to exactly two of them, and so $b_1$ and $b_2$ have no
common neighbour in $K$. The same holds for the
pairs from $A$ and $C$. This shows that the two vertices adjacent to $b_1$ cannot both lie in $A$,
for otherwise $a_1$ and $a_2$ would be both adjacent to $b_1$. By symmetry
each vertex is adjacent to exactly one element of the other non-singleton kernel classes. Each vertex with its two neighbours form a transversal for $\{A,B,C\}$; thus there are only
 two possible  type of configurations: either the edges form two disjoint $3$-cycles,
both of which intersect all of $A$, $B$, $C$ or the edges form a $6$-cycle that transverses $A$, $B$, $C$ in a periodic order. In different words, we can assume without loss of generality that we have $a_1-b_1-c_1$ and $a_2-b_2-c_2$. Thus, either $a_1-c_1$ (and hence we have two $3$-cycles $a_1-b_1-c_1-a_1$ and $a_2-b_2-c_2-a_2$, see Figure \ref{fi:3transp1}), or $a_1-c_2$ and we have one $6$-cycle $a_1-b_1-c_1-a_2-b_2-c_2-a_1$ (see Figure \ref{fi:3transp2}).

\begin{figure}[h]
\[
\xy
(-5,5)*{}="b";
(5,5)*{}="c";
(-5,25)*{}="f";
(5,25)*{}="g";
"b";"f" **\dir{-};
"b";"c" **\dir{-};
"c";"g" **\dir{-};
"g";"f" **\dir{-};
(15,5)*{}="bS";
(25,5)*{}="cS";
(15,25)*{}="fS";
(25,25)*{}="gS";
"bS";"fS" **\dir{-};
"bS";"cS" **\dir{-};
"cS";"gS" **\dir{-};
"gS";"fS" **\dir{-};
(35,5)*{}="bT";
(45,5)*{}="cT";
(35,25)*{}="fT";
(45,25)*{}="gT";
"bT";"fT" **\dir{-};
"bT";"cT" **\dir{-};
"cT";"gT" **\dir{-};
"gT";"fT" **\dir{-};
(0,10)*{\bullet}="d";
(0,0)*{}="a";
(40,10)*{\bullet}="d2";
(40,0)*{}="a2";
(20,10)*{\bullet}="d1";
"a";"a2" **\dir{-};
"a2";"d2" **\dir{-};
"d2";"d" **\dir{-};
"d";"a" **\dir{-};
%
(0,30)*{}="dT";
(0,20)*{\bullet}="aT";
(20,20)*{\bullet}="a1T";
(40,30)*{}="d2T";
(40,20)*{\bullet}="a2T";
(20,30)*{}="d1T";
"aT";"a2T" **\dir{-};
"a2T";"d2T" **\dir{-};
"d2T";"dT" **\dir{-};
"dT";"aT" **\dir{-};
\endxy
\]
\caption{One of the two possible induced subgraphs on $K$}\label{fi:3transp1}
\end{figure}
\begin{figure}[h]
\[
\xy
(-5,5)*{}="b";
(5,5)*{}="c";
(-5,25)*{}="f";
(5,25)*{}="g";
"b";"f" **\dir{-};
"b";"c" **\dir{-};
"c";"g" **\dir{-};
"g";"f" **\dir{-};
(15,5)*{}="bS";
(25,5)*{}="cS";
(15,25)*{}="fS";
(25,25)*{}="gS";
"bS";"fS" **\dir{-};
"bS";"cS" **\dir{-};
"cS";"gS" **\dir{-};
"gS";"fS" **\dir{-};
(35,5)*{}="bT";
(45,5)*{}="cT";
(35,25)*{}="fT";
(45,25)*{}="gT";
"bT";"fT" **\dir{-};
"bT";"cT" **\dir{-};
"cT";"gT" **\dir{-};
"gT";"fT" **\dir{-};
(0,20)*{\bullet}="e";
(0,10)*{\bullet}="d";
(0,0)*{}="a";
(40,10)*{\bullet}="d2";
(40,0)*{}="a2";
(20,10)*{\bullet}="d1";
(20,20)*{\bullet}="a1T";
(40,20)*{\bullet}="a2T";
(0,30)*{}="dT";
(40,30)*{}="d2T";
"a";"a2" **\dir{-};
"a2";"d2" **\dir{-};
"d2";"d1" **\dir{-};
"d1";"e" **\dir{-};
"e";"dT" **\dir{-};
"dT";"d2T" **\dir{-};
"d";"a1T" **\dir{-};
"a1T";"a2T" **\dir{-};
"d2T";"a2T" **\dir{-};
"a";"d" **\dir{-};
\endxy
\]
\caption{One of the two possible induced subgraphs on $K$}\label{fi:3transp2}
\end{figure}

In either case,  the triple transposition $g=(a_1 \, a_2)(b_1 \, b_2)(c_1 \, c_2)$
is an automorphism of the induced subgraph on $K$. In fact, since $N(a_1) \cap \overline{K} = N(a_2) \cap\overline{K}$,  $N(b_1) \cap \overline{K} = N(b_2)\cap \overline{K}$,  $N(c_1) \cap \overline{K} = N(c_2) \cap \overline{K}$,
the trivial extension of $g$ is an automorphism of $\Gam$.
As explained in Subsection~\ref{smallmindegree}, this is impossible. So we have:

\begin{theorem}
Let $G$ act primitively on $X$, and let $f \in T(X)$ have kernel type  $(2,2,2,1,\dots,1)$.
Then $G$ synchronizes $f$.
\end{theorem}

With the results above about transformations of
kernel type $(3,2,1,\dots,1)$ (taking $p=3$ in Theorem \ref{kerp211}) and the results from \cite{arcameron22}
 ($k=4$ in Theorem 2) about transformations of kernel type $(4,1,\dots,1)$,
we get Theorem \ref{main}.

\section{Sets with a small neighbourhood}\label{k-2sets}


In this section, we will exploit sets of vertices that share large number of adjacent vertices to show that certain kernel types are always synchronized. A consequence of the results in this section is that primitive groups synchronize every transformation of kernel type $(p,3,1, \dots,1)$, for $p\ge 4$. We also introduce some notation that will be very important in the next section.

Let $\Gam$ be a regular graph with valency $k$, and let $A \subseteq V(\Gam)$. We
 say that $A$ is a \emph{small neighbourhood set of defect} $d$  if $|\cup_{a \in A} N(a)| \le k+d$.

We will assume throughout this section that $\Gam$ is a graph with primitive automorphism group and clique number equal to its chromatic number.

\begin{lemma} \label{lm:ANA} Assume that $A$ is a small neighbourhood set of defect $2$ in $\Gam$ of size $l \ge 3$. Set
$N_A=\cup_{a \in A} N(a)$. Let $x,y,z \in A$ be distinct, and $w \in N_A$. Then
\begin{enumerate}
\item $|N(x) \cap N(y)|=k-2$,
\item  $N(x) \cup N(y)=N_A$,
\item $N(z) \subseteq N(x) \cup N(y)$,
\item $|N(w) \cap A| \ge l-1$, \label{NA}
\item the $2$ elements of $N_A \setminus N(z)$ are in $N(x) \cap N(y)$,  \label{extra2}
\item $N_A$ contains at least $2 l$ elements that are not adjacent to all elements of $A$.
\end{enumerate}
\end{lemma}
\begin{pf}
The first two claims follow from $|N(x) \cup N(y)| \le k+2$ in connection with the pigeonhole principle. The third follows from the second. For the fourth, assume that $x,y \in A \setminus N(w)$, $x \ne y$. Then
$N(x) \cup N(y) \ne N_A$, as $w \notin N(x) \cup N(y)$, for a contradiction. For (\ref{extra2}) notice that
any counterexample $w$ would contradict (\ref{NA}). The last claim now follows from (\ref{extra2}).
\qed
\end{pf}

Define $l_1=2$ and $l_d= l_{d-1}+d$ for $d \ge 2$.
\begin{lemma}\label{lm:smnbound}
For $d \ge 1$, $\Gam $ does not contain any small neighbourhood set $A$ of defect $d$ and size $l_d$.
\end{lemma}
\begin{pf}
 The proof is by induction on $d$. For $d=1$, notice that a small neighborhood set $A$ of defect $1$ and size $2$ contradicts Corollary \ref{cork-2} in connection with the pigeonhole principle.

 So let $d \ge 2$ and assume that the result holds for smaller values of $d$. By way of contradiction let $A$ be a small neighbourhood set
 of defect $d$ with $l_d$ distinct elements. We may assume that $|N_A|= k+d$. Let $w \in N_A$. We claim that $|N(w) \cap A|> l_d- l_{d-1}$. Indeed, if $A' \subseteq A \setminus  N(w)$ for some $A'$ with $|A'|=l_{d-1}$, then
 $\cup_{a \in A'}N(a) \subseteq N_A \setminus\{w\}$, and $A'$ would be a small neighbourhood set of defect at most $d-1$ and size $l_{d-1}$. Such an $A'$ does not exist by our inductive assumption, and so
 $|N(w) \cap A|> l_d- l_{d-1}$.

 Now let $x,y \in A$, and $g \in G$ such that $xg \in N_A$, $yg \notin N_A$. Such $g$ exists by primitivity. As $N(x) \cup N(y)$ is contained in a set of size $k+d$, $N(y)$ intersects any subset of $N(x)$ of size $d+1$.  The same holds for $N(yg)$ and $N(xg)$. As $xg \in N_A$, there are at least $l_d-l_{d-1}+1 =d+1$ elements in $A \cap N(xg)$. Hence $yg$ is connected to one of those elements, and hence in $N_A$, for a contradiction.
 \qed
\end{pf}

We have everything needed to prove the main theorem of this section.

\begin{theorem}\label{t:dominantclass} Let $d_1, \dots, d_j\ge 2$ be integers and $d=-j+\Sigma d_i$. Let $l \ge l_d$. Then $G$ synchronizes every map $f$ of kernel type $(l,d_1,d_2,\dots,d_j,1, \dots, 1)$.
\end{theorem}
\begin{pf} Assume otherwise, and let  $A$ be the kernel class of size $l$ of $f$, $B_i$ be the other non-singleton kernel classes of $f$, and $x \in A$. In $\Gam$, the
 elements of $N_A$ all map to $N(xf)$ of size $k$. This set has at most $k+(d_1-1)+(d_2-1) + \dots + (d_j-1) =k+d$ preimages. It follows that $A$ is a small neighbourhood set of defect $d$ and size at least $l_d$, contradicting Lemma \ref{lm:smnbound}.
 \qed
\end{pf}
Theorem \ref{t:dominantclass} is applicable if $j=1$ and $d=3$, in which case $l_d=4.$ In Subsection \ref{ker3311}, we will show that a primitive permutation group synchronizes every map of kernel type $(3,3,1,\dots,1)$. Together,
these results imply the following corollary (see also Theorem \ref{kerp211}).
\begin{cor} Let $p \ge 3$, and $G$ a primitive permutation group on $X$. Then $G$ synchronizes every transformation on $X$ of kernel type $(p,3,1, \dots,1)$.
\end{cor}

\section{Maps of rank $n-4$}\label{largernk2}

The aim of this section is to prove the following:

\begin{theorem}\label{th:n-4}
Let $G$ be a primitive group acting on a set of vertices $X$ with $|X|=n \ge 5$. Then $G$ synchronizes every map of rank $n-4$.
\end{theorem}
We will first  prove various auxiliary lemmas and describe our general proof strategy.
The actual proofs involve a large number of subcases and will be divided over the next three subsections, each of
which covers a particular kernel class.

Throughout our proof of Theorem \ref{th:n-4}, we assume that $G$ is a
primitive group of degree $n$ over a set $X$, $f$ is a transformation of
rank  $n-4$, and $G$ does not synchronize $f$.  We let $\Gam'=\gr'(S)$
be the graph constructed earlier for $S=\langle G\cup f \rangle)$, $k$ be the
valency of $\Gam'$, and $r$ its clique size.

The five possible kernel classes for a map of rank $n-4$ are $(5,1, \dots, 1)$, $(4,2,1,\dots,1), (3,3,1, \dots,1), (3,2,2,1,\dots,1)$, and $(2,2,2,2,1,\dots,1)$. If $f$ is  one of the first two types, the result was shown in
\cite{arcameron22} and Theorem \ref{kerp211}. The remaining three cases are covered in the following subsections.

For each kernel type, we will denote by $K$ the union of the non-singleton kernel classes of $f$. For any given
non-singleton kernel class $Z$, we let $N_Z=\cup_{z \in Z} N(z)$, and let $N_Z'= N_Z \cap \overline{K}$. We repeat that for all such $Z$, $|N_Z|\ge k+2$, as  neighbourhoods of distinct elements in $Z$ may only have intersection of size at most $k-2$.

We will distinguish several cases by the induced subgraph on the set $Kf$. Let $Z$ be a non-singleton
kernel class of $f$ with image $z'$. Let $Y_1, \dots, Y_m$ be those non-singleton kernel classes that
map to neighbours of $z'$. We refer to the number $p_Z=|Y_1|+\dots + |Y_m|+(k-m)$ as the \emph{number
of potential neighbours} of $Z$, and to $p_Z'=k-m$ as the \emph{number of potential singleton kernel class neighbours} of $Z$.

\begin{lemma}
$|N_Z| \le p_Z$, $|N_Z'|\le p_Z'$.
\end{lemma}
\begin{pf}
$z'$ has $m$ neighbours that are images of non-singleton kernel classes and hence $k-m$ neighbours that are either images of singleton kernel classes or not in the image of $f$. If $z \in Z$ and $y$ is such that $z-y$,
then $y$ must be a preimage of a neighbour of $z'$, hence  $y \in Y_i$ for some $i$ or $y$ is the singleton class preimage of one
of remaining $k-m$ elements of $N(z')$. The results follow.\qed
\end{pf}
For $z \in X$, let $[z]$ denote the kernel class of $f$ containing $z$.
\begin{lemma}\label{lm:estimatelow}
Let $r$ be the number of edges in the induced subgraph of $K$. Then
$r\ge \frac{1}{2}\Sigma_{z \in K} (k-p_{[z]}')$.
\end{lemma}
\begin{pf}
For any given $z \in K$, all neighbours of $z$ that lie in singleton kernel classes are in $N_{[z]}'$. By the previous
lemma $|N_{[z]}'| \le p_{[z]}'$. Hence $z$ has at least $k-p_{[z]}'$ neighbours in $Z$. Summing over all
$z \in K$,
we obtain a lower bound on the number of pairs in the adjacency relation on $K$. The result follows.\qed
\end{pf}
\begin{lemma}\label{lm:preestimate}
Suppose that there are $s$ non-singleton kernel classes, and that the induced subgraph on $Kf$ has $r'$ edges.
Let $r$ be the  number of edges in $K$, then
  $$r \le sk-r'-\Sigma |N_Z'|+6,$$
where the sum is over the non-singleton kernel classes $Z$ of $f$.
\end{lemma}
\begin{pf} Consider the two induced graphs on $X$ and  $Xf$. We will estimate the difference
in their number of edges in two ways.

$Xf$ is obtained from $X$ by deleting $4 $ vertices, namely the non-images of $f$. Each of
these is a vertex of $k$ edges. Hence we lose $4k$ edges minus the number that we count twice because
both of their vertices are non-images of $f$. There are at most $6$ such edges between $4$ vertices.
 Hence we lose at least $4k-6$ edges.

We obtain another estimate by comparing various subsets of edges and their images under $f$. We start with
those edges that are within $K$: here $r$ edges are mapped onto $r'$ edges for a loss
of $r-r'$.

For each non-singleton kernel class, $Z$ let $r_Z$ be the number of edges between $Z$ and $K\setminus Z$.
Then there are $|Z|k-r_Z$ edges between $Z$ and elements in singleton kernel classes.
These edges map to the $|N_Z'|$ edges between the image of $Z$ and the images of $N_Z'$. Hence
we have an effective loss of $|Z|k-r_Z -|N_Z'|$ of edges.

Finally we note that all edges between singleton classes are mapped injectively to other edges, so we do not
encounter any loss for them.

Summing up, we obtain a loss of at most
\begin{eqnarray*}
  (r-r') + \Sigma \left( |Z|k-r_Z-|N_Z'|\right) &=& r-r' +\Sigma|Z|k - \Sigma r_Z -\Sigma |N_Z'| \\
   &=& r-r'+|K|k - 2r-\Sigma |N_Z'| \\
   &=& |K|k -r-r'- \Sigma|N_Z'|
\end{eqnarray*}
edges, where the sums are over the set of non-singleton kernel classes indexed by $Z$. Comparing with the
lower bound $4k-6$, we get that
\begin{eqnarray*}
   r&\le& (|K|k-r'-\Sigma |N_Z'|)-(4k-6) \\
   &=& (|K|-4)k -r' -\Sigma |N_Z'|+6 \\
  &=& sk-r'-\Sigma |N_Z'|+6.
\end{eqnarray*}
\qed
\end{pf}
In the following, we will only be dealing with kernel classes $Z$ that satisfy $p_Z \in \{k+2, k+3\}$. As
$p_Z \ge |N_Z| \ge k+2$, in cases where $p_Z=k+2$,  we get that $p_Z=|N_Z|$. Hence
every potential neighbour of $Z$ is in fact a neighbour. In particular, every potential singleton class neighbour
is also a neighbour, which implies that
$|N_Z'|=p_Z'=k-m_Z$, where $m_Z$ is the number of neighbours of the image of $Z$ in $Kf$.
In case that $p_Z=k+3$, one potential neighbour might not be a neighbour (or might not exist, if the image of $Z$ has a neighbour that
is not in the image of $f$). Hence in this case $|N_Z'| \in \{p_Z',p_Z'-1\}$.

\begin{lemma}\label{lm:estimate}
Under the conditions of Lemma \ref{lm:preestimate}, assume that for all non-single\-ton kernel classes $Z$,
$p_Z \in \{k+2, k+3\}$. Let
$d$ be the number of kernel classes for which $p_Z=k+3$.
Then $r \le r' +d+6$.

Moreover, for each $Z$, let $m_Z$ be the number of neighbours of the image of $Z$ that lie in $K f$.
If $r=r' +i+6$, for some $1\le i \le d$, then there are at least $i$ non-singleton kernel classes $Z$ for which $|N_Z'| = p_Z'-1=k-m_Z -1$.
\end{lemma}
\begin{pf}
By Lemma \ref{lm:preestimate}, $r \le sk-r'-\Sigma |N_Z'|+6$, and as pointed out after the lemma, we have $|N_Z'|=k-m_Z$,
if $p_z=k+2$, or $|N_Z'| \ge k-m_Z-1$, if $p_z=k+3$. Assume that there are exactly $j$ kernel classes $Z$ for which $|N_Z'| = k-m_Z-1$. Then
\begin{eqnarray*}
  r &\le& sk-r'-\Sigma |N_Z'|+6 \\
   &= & sk-r' -\left( \Sigma (k-m_Z) -j\right)+6 \\
   &=& sk-r'-\left(sk-2r'\right) +j+6\\
   &=& r'+j+6.
\end{eqnarray*}
$|N_Z'| =  k-m_Z-1$ implies that $p_Z=k+3$, therefore $j \le d$, and the first statement of the lemma follows. Assuming $r=r'+i+6$, we obtain $i \le j$, which shows the second statement.
\qed
\end{pf}
Our proof of Theorem \ref{th:n-4} proceeds by considering for each kernel class all potential combinations of induced subgraphs on $K$ and $Kf$.
All configurations whose number of edges lie within the bounds of Lemmas \ref{lm:estimatelow} and \ref{lm:estimate} will be further restricted
and eventually excluded.

One of our most common arguments will be to construct  a contradiction to Lemma
\ref{primgr}. As we will use this construction extensively, we will introduce some special notation for it.
By a CME -- standing for \emph{clique minus one edge} --  we mean a set of vertices of size $r+1$ that contains at most one non-edge, i.e., a configuration
that violates either Lemma \ref{primgr} or the fact that $r$ is the clique number of $\Gam'$.

For distinct  vertices $x,y,z$, with $x-y$, the expression CME$(x-y, z)$ means that for any $r$-clique $L$ that contains
the edge from $x$ to $y$ (whose existence follows from the definition of $\Gam'$), the set $L \cup \{z\}$ is
a CME. A typical application will be that
 $z$ is in the same kernel class as one of $x$ or $y$, and adjacent to the other one. Often
we will have that $N_{[z]}' \subseteq N(z)$ due to
$z$ having not enough neighbours in $K$ to omit a vertex from $N_{[z]}'$. It then just remains to check
that all vertices in $K$ adjacent to both $x$ and $y$ are also adjacent to $z$.

Another tool is to utilize small neighbourhood sets of defect $2$. We always have such a set of size at least $2$ available if
we have a kernel class $Z$ with $p_Z=k+2$. By transitivity of $G$, every element is then part of such a set.  The following
lemmas draw consequences in these cases.
\begin{lemma} \label{lm:defect2b} Let $x,y \in X$, such that in $\Gam$, $|N(x) \cap N(y)| =k-2$.
\begin{enumerate}
\item $x$ and $y$ are non-adjacent.\label{nonadj}
\item Suppose that $Z$ is a kernel class of $f$ such that $N(x) \cap Z \ne \emptyset \ne N(y) \cap Z$, but that $Z \cap N(x) \ne Z \cap N(y)$. Then $xf =yf$.  \label{mapclique}
\end{enumerate}
\end{lemma}
\begin{pf}
Assume that $x$ and $y$ are adjacent. Then $x \in N(y) \setminus N(x)$. Let $z$ be the other element of $N(y) \setminus N(x)$.
By Proposition \ref{pr:derived}, there exists
an $r$-clique $L$ containing $y$, but not containing $x$.
 Then $L \cup \{x\}$ is a CME, as it has $r+1$ elements and at most one non-edge between $x$ and $z$.
By contradiction, we obtain (\ref{nonadj}).

Now  in the situation of (\ref{mapclique}), say w.l.o.g. that $z \in (Z\cap N(y))\setminus N(x)$. Let $L$ be an $r$-clique containing $y$ and avoiding the unique element in $N(y)\setminus \left(N(x)\cup \{z\}\right)$. Then $z \in L$ for otherwise $L'=L \cup\{x\}$ is a CME. Hence $z \in L'$, and $L'$ is missing two edges, namely $(x,y)$ and $(x,z)$. Now
$Z \cap N(x) \ne \emptyset$, hence there is an edge from $x$ to an element of $Z$, and hence the non-edge $(x,z)$ maps to the edge $(xf, zf)$. It follows that $L'f$ cannot have
$r+1$ elements, for otherwise it would be a CME. So $f$ must identify two elements of $L'$. These cannot be any elements of the clique $L$. $x$ is adjacent to all elements of $L \setminus \{y,z\}$,
and $(xf,zf)$ is an edge. Thus $xf=yf$ by elimination.\qed
\end{pf}
\begin{lemma}\label{lm:smnd2s2pre}
Suppose that in $\Gam'$ we have a small neighbourhood set of defect $2$ and size at least $2$. Then there exist vertices $x,y,z \in \Gam'$ such that
 $|N(x) \cap N(y)|=k-2$, $|N(y) \cap N(z)|=k-2$, $|N(x) \cap N(y) \cap N(z)| <k-2$. Moreover, such triples exist for any chosen vertex $y$.
\end{lemma}
\begin{pf}
Let $\sim$ be the relation on $\Gam'$ defined by $x \sim y$ if either $x=y$ or $|N(x) \cap N(y)|=k-2$. The relation $\sim$ is clearly reflexive, symmetric, and preserved by  $G$.

Assume that for all $x,y,z \in \Gam'$, $|N(x) \cap N(y)|=k-2=|N(y) \cap N(z)|$ implies that
 $|N(x) \cap N(z)|=k-2$. Our assumption means that
$ \sim $ is transitive and hence a $G$-compatible
equivalence relation on $X$. By primitivity of $G$, $\sim$ is trivial or universal. However, $\sim $ is non-trivial as we assumed that $\Gam'$ has a small neighbourhood set of defect $2$,
and it is not universal, as adjacent elements of $\Gam'$ are not in $\sim$ by Lemma \ref{lm:defect2b}(\ref{nonadj}). By contradiction, there
exist $x,y,z \in \Gam'$, with $|N(x) \cap N(y)|=k-2=|N(y) \cap N(z)|$, and $k-2 > |N(x) \cap N(z)| \ge |N(x) \cap N(y) \cap N(z)|$.

The last assertion follows from the transitivity of $G$. \qed
\end{pf}
\begin{lemma}\label{lm:smnd2s2} Suppose that in $\Gam'$ we have a small neighbourhood set of defect $2$ and size at least $2$. Let $y \in \Gam'$, and $y', \bar y\in N(y), y' \ne \bar y$. Then there exists a $w \in \Gam'$
such that $|N(y) \cap N(w)| =k-2$ and $N(w) \cap \{y', \bar y\} \ne \emptyset$.
\end{lemma}
\begin{pf}
Given $y$, let $x,z$ be the elements constructed in Lemma \ref{lm:smnd2s2pre}. Then $|N(x) \cap N(y)|=|N(z) \cap N(y)|=k-2$. We claim that  one of $x,z$ is adjacent to an element of $\{y', \bar y\}$.
For assume otherwise, then $|N(x) \cap N(y) \cap N(z)|= |N(y) \setminus \{y',\bar y\}|=k-2$, contradicting Lemma  \ref{lm:smnd2s2pre}. The result follows. \qed
\end{pf}
\begin{lemma}\label{lm:Nimage}
Let $x,y \in \Gam'$, $xf \ne yf$, such that $\{x,y\}$ is a small neighbourhood set of defect 2. Let $N= N(x) \cap N(y).$ If for every non-singleton kernel class $Z$ of $f$,
$|N\cap Z| \le 1$, then $xf$ and $yf $ are non-adjacent.
\end{lemma}
\begin{pf} As $\{x,y\}$ is a small neighbourhood set of defect 2, $|N|=k-2$. Consider $Nf$. As $|N \cap Z| \le 1$ for all kernel classes $Z$, $f$ maps $N$ injectively, and so $|Nf|=k-2$. Moreover, $x,y$ are adjacent to every
element in $N$, and as $xf \ne yf$, $N\cup \{x,y\}$ is mapped injectively by $f$, as well. It follows that $|N(xf) \cap N(yf)| \ge |Nf|=k-2$, which implies that $\{xf,yf\}$ are also a small neighbourhood set of defect 2.
The result now follows with Lemma \ref{lm:defect2b} (\ref{nonadj}).
\qed
\end{pf}

\begin{lemma}\label{lm:defect2c}
Let $A_1$, $A_2$ be small neighbourhood sets of defect $2$ and size $3$. If $|A_1 \cap A_2| \ge 2$ then $A_1=A_2$.
\end{lemma}
\begin{pf}
Let $A_1=\{x,y,z_1\}, A_1=\{x,y,z_2\}$. Then
$$k+2\le|N_{\{x,y\}}| \le |N_{A_1}|=k+2,$$
and so $N_{\{x,y\}}= N_{A_1}$. Symmetrically,
$N_{A_2}=N_{\{x,y\}}=N_{A_1}$ which implies that $N_{A_1 \cup A_2}=N_{A_1}$, and so $|N_{A_1 \cup A_2}|=k+2$.
By Lemma \ref{lm:smnbound}, there are no small neighbourhood sets of
defect $2$ and size $4$, hence $z_1=z_2$ and $A_1=A_2$.\qed
\end{pf}

\subsection{Maps of kernel type $(3,3,1, \dots,1)$}\label{ker3311}

Let $f$ be a map of kernel type $(3,3,1, \dots,1)$, $A=\{a_1,a_2,a_3\}$, $B =\{b_1,b_2,b_3\}$ the non-singleton kernel classes of $f$, and
assume that $\Gam$ has $r$ edges between $A$ and $B$. In order for $p_A \ge  k+2$, the images of
$A$ and $B$ need to be connected and we get $k+2=p_A=|N_A|=p_B=|N_B|$
and $|N_A'|=|N_B'|=k-1$.
Hence $A, B$ are small neighbourhood sets of defect $2$.

Our next goal is to bound $r$. By Lemma \ref{lm:ANA}, every element of $N_A$ is
 adjacent to at least $2$ elements in $A$, hence $r \ge 6$. Lemma \ref{lm:estimate} shows that $r\le 7$.

\begin{figure}[h]
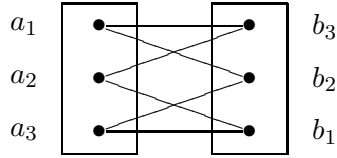

\[
\xy
(-5,5)*{}="b";
(5,5)*{}="c";
(-5,25)*{}="f";
(5,25)*{}="g";
"b";"f" **\dir{-};
"b";"c" **\dir{-};
"c";"g" **\dir{-};
"g";"f" **\dir{-};
(15,5)*{}="bS";
(25,5)*{}="cS";
(15,25)*{}="fS";
(25,25)*{}="gS";
"bS";"fS" **\dir{-};
"bS";"cS" **\dir{-};
"cS";"gS" **\dir{-};
"gS";"fS" **\dir{-};
(0,22)*{\bullet}="a";
(0,8)*{\bullet}="e";
(0,15)*{\bullet}="d";
(20,22)*{\bullet}="a1";
(20,8)*{\bullet}="e1";
(20,15)*{\bullet}="d1";
"a";"a1" **\dir{-};
"a";"d1" **\dir{-};
"a1";"d" **\dir{-};
"d";"e1" **\dir{-};
"e";"d1" **\dir{-};
"e";"e1" **\dir{-};
(-10,22)*{a_1}="a";
(-10,8)*{a_3}="e";
(-10,15)*{a_2}="d";
(30,22)*{b_3}="a1";
(30,8)*{b_1}="e1";
(30,15)*{b_2}="d1";
\endxy
\]
\caption{The induced subgraph on $K$ with $6$ edges}\label{fi:ker33ed6}
\end{figure}

\begin{figure}[h]
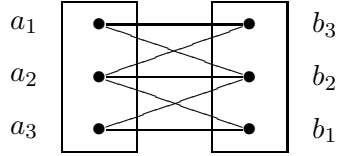

\[
\xy
(-5,5)*{}="b";
(5,5)*{}="c";
(-5,25)*{}="f";
(5,25)*{}="g";
"b";"f" **\dir{-};
"b";"c" **\dir{-};
"c";"g" **\dir{-};
"g";"f" **\dir{-};
(15,5)*{}="bS";
(25,5)*{}="cS";
(15,25)*{}="fS";
(25,25)*{}="gS";
"bS";"fS" **\dir{-};
"bS";"cS" **\dir{-};
"cS";"gS" **\dir{-};
"gS";"fS" **\dir{-};
(0,22)*{\bullet}="a";
(0,8)*{\bullet}="e";
(0,15)*{\bullet}="d";
(20,22)*{\bullet}="a1";
(20,8)*{\bullet}="e1";
(20,15)*{\bullet}="d1";
"a";"a1" **\dir{-};
"a";"d1" **\dir{-};
"a1";"d" **\dir{-};
"d";"e1" **\dir{-};
"d";"d1" **\dir{-};
"e";"d1" **\dir{-};
"e";"e1" **\dir{-};
(-10,22)*{a_1}="a";
(-10,8)*{a_3}="e";
(-10,15)*{a_2}="d";
(30,22)*{b_3}="a1";
(30,8)*{b_1}="e1";
(30,15)*{b_2}="d1";
\endxy
\]
\caption{The induced subgraph on $K$ with $7$ edges}\label{fi:ker33ed7}
\end{figure}

We will treat the two cases $r=6,7$ simultaneously. If  $r=6$, then every element of $B$ is adjacent to exactly $2$ elements of $A$ and vice versa. If $r=7$ then  exactly one element of $A$ is adjacent to all vertices in $B$,
exactly one element of $B$ is adjacent to all vertices in $A$, and the
remaining elements of $A \cup B$ have exactly $2$ neighbours in $K$. Hence w.l.o.g., we may assume
 that all edges in $A \cup B$ lie on the $6$-cycle $a_1-b_3-a_2-b_1-a_3-b_2-a_1$, except for potentially an extra edge between $a_2$ and $b_2$ in case that $r=7$. These two configurations are depicted in Figures \ref{fi:ker33ed6} and
 \ref{fi:ker33ed7}.

\begin{figure}[h]
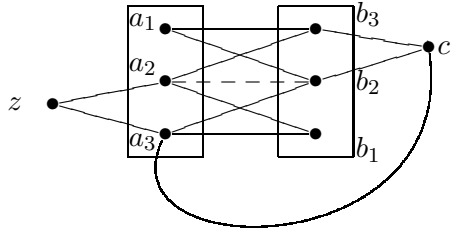

\[
\xy
(-5,5)*{}="b";
(5,5)*{}="c";
(-5,25)*{}="f";
(5,25)*{}="g";
"b";"f" **\dir{-};
"b";"c" **\dir{-};
"c";"g" **\dir{-};
"g";"f" **\dir{-};
(15,5)*{}="bS";
(25,5)*{}="cS";
(15,25)*{}="fS";
(25,25)*{}="gS";
"bS";"fS" **\dir{-};
"bS";"cS" **\dir{-};
"cS";"gS" **\dir{-};
"gS";"fS" **\dir{-};
(0,22)*{\bullet}="a";
(0,8)*{\bullet}="e";
(-15,12)*{\bullet}="Z";
(-20,12)*{z};
(0,15)*{\bullet}="d";
(20,22)*{\bullet}="a1";
(20,8)*{\bullet}="e1";
(20,15)*{\bullet}="d1";
"Z";"e" **\dir{-};
"Z";"d" **\dir{-};
"a";"a1" **\dir{-};
"a";"d1" **\dir{-};
"a1";"d" **\dir{-};
"d";"e1" **\dir{-};
"d";"d1" **\dir{--};
"e";"d1" **\dir{-};
"e";"e1" **\dir{-};
(-3,23)*{a_1};
(-3,7)*{a_3};
(-3,17)*{a_2};
(27,24)*{b_3};
(27,6)*{b_1};
(27,15)*{b_2};
(37,19.5)*{c};
(35,19.5)*{\bullet}="C";
"C";"d1" **\dir{-};
"C";"a1" **\dir{-};
"e";"C" **\crv{(-10,-10) & (40,-10)} ?>* \dir{-};

\endxy
\]
\caption{The construction from Lemma \ref{lm:331struc}.\label{fi:ker33lemma}}
\end{figure}

\begin{lemma} \label{lm:331struc}
There exist  unique elements $z \in N_A', c \in N_B'$ that are not adjacent to $a_3, b_1$, respectively.  Moreover, $c$ is adjacent to $a_3$.

\end{lemma}
\begin{pf}
We have that $|N(b_1) \cap A|=2$. It follows that $|N(b_1) \cap N_B'|=k-2$. As $|N_B'|=k-1$, there is exactly one element $c$  in $N_B'$ that is not connected to $b_1$. The existence and uniqueness of $z$ follow symmetrically. By (\ref{NA}) of Lemma \ref{lm:ANA}, we have the edges $b_3-c-b_2$, and $a_2-z-a_3$.

Now consider an $r$-clique $L$ containing the edge $a_3-b_2$. We have that $b_1-a_3$ and
$L \setminus \{a_3,b_2\}\subseteq N_B'\subseteq N(b_1) \cup\{c\}$.
It follows that $c \in L$ for otherwise $L \cup \{b_1\}$ would be a CME, missing only an edge between $b_1$ and $b_2$. Hence $c-a_3$.

The construction from this lemma is depicted in Figure~\ref{fi:ker33lemma}.
Note that there may be additional edges that are not depicted, except for the confirmed non-edges $(c,b_1)$, $(z,a_1)$.
The dotted edge is the additional edge in the case $r=7$.\qed
\end{pf}

%
%
%

Let $g\in G$ be such that $a_1g \in A$, $a_3g \notin A$. Consider $A'=Ag^{-1}$. It is a small neighbourhood set of defect $2$, as $A$ has this property. Moreover $a_1 \in A \cap A'$ but $ A \ne A'$, as $a_3 g \notin A$. By Lemma \ref{lm:defect2c}, $A' \cap A=\{a_1\}$. Let $A' =\{a_1,x,y\}$. $b_3 \in N(a_1)$ and hence by (\ref{NA}) of Lemma \ref{lm:ANA}, one element of $x,y$, say $x$, must be adjacent to $b_3$. Hence $x \in N_B\setminus A$.

As $xf \ne a_1f$, by Lemma \ref{lm:defect2b}, $N(a_1) \cap B= N(x) \cap B=\{b_2,b_3\}$. As $b_1 \notin N(x)$, $x=c$, where $c$ is from Lemma \ref{lm:331struc}. By the same lemma, we have that $x-a_3$. As $x$ is not adjacent to $a_1$, we have that $x=c=z$ once again by Lemma \ref{lm:331struc}.

Now, consider the third element $y$ of $A'$. If $y$ would be adjacent to $b_2$ or $b_3$, repeating the argument from the previous paragraph yields $y=c$. As $x \ne y$, it follows that $y$ is not adjacent to
 $b_2$ or $b_3$. As $$|N(c) \cup N(y)|=|N(x) \cup N(y)|=k+2,$$
 it follows that $y$ is adjacent to every element in $N(c) \setminus \{b_2,b_3\}$. So $y \in N(a_3)$ and hence $y \in N_A$. As $y \notin N(a_1)$ the uniqueness of $z=x$ implies that $y \notin N_A'$.
So $y \in B$, and hence $y=b_1$, as $b_2, b_3$ are adjacent to $a_1$. It follows that $\{a_1,b_1\}$ is a small neighbourhood set of defect $2$.

Consider $N=N(a_1) \cap N(b_1)$ of size $k-2$. $N$ has no elements in $K$, and hence $|N \cap Z|\le 1$ for all kernel classes $Z$ of $f$. By Lemma \ref{lm:Nimage}, $a_1 f$ and $b_1 f$ are non-adjacent; however, this
 is false in our construction.


Our assumption was that $G$ synchronizes the transformation $f$. Hence by contradiction, Theorem \ref{th:n-4} holds for transformations of kernel type $(3,3,1, \dots,1)$.

\subsection{Transformations of kernel type $(3,2,2,1, \dots,1)$}
Let $A=\{a_1,a_2\}$, $B =\{b_1,b_2,b_3\}$, $C=\{c_1,c_2\}$ be the non-singleton kernel classes of $f$.
The requirement that $p_Z \ge k+2$ for all kernel classes $Z$ implies that $K f$ is connected. Hence the induced graph on $Kf$ is a $2$-path or a triangle.

\paragraph{The induced graph on $Kf$ is a $2$-path} $\phantom{.}$

The requirement that $p_B \ge k+2$ implies that there must be edges from $B$ to both $A$ and $C$, hence $a_1f - b_1f -c_1f$.

Let $r$ be the number of edges in $K$, by Lemma \ref{lm:estimate} we conclude that $r \le 8$.
As $|N_B'|=k-2$, each element of $B$ has at least $2$ neighbours in $K$. Together, these constraints imply that at least one element of $B$ has exactly $2$ neighbours in $K$. If this holds for all elements of $B$,
then for at least two distinct $x,y \in B$,  $N(x) \cap N(y)\cap K \ne \emptyset$. Otherwise, there are $x, y \in B$, with $|N(x) \cap K|=2$, $|N(y) \cap K|\ge 3.$
In both cases,  $x, y \in B$ satisfy $|N(x) \cap K|=2$, and $N(x) \cap N(y)\cap K \ne \emptyset$. Say w.l.o.g. that $x= b_1$, $y=b_2$, and $b_1-c_1-b_2$.

 We claim that we have a CME$(b_2-c_1,b_1)$.
For let $L$ be an $r$-clique containing $b_2,c_1$, then  $L\setminus \{b_2,c_1\} \subseteq N_B'$. However $N_B' \subseteq N(b_1)$, as $|N_B'|=k-2$, and $b_1$ has only two neighbours in $K$. Hence $L \cup \{b_1\}$
is only missing one edge between $b_1$ and $b_2$, and is a CME.

By contradiction, we can exclude the case that $Kf$ is a $2$-path.

\paragraph{The images of the non-trivial kernel classes form a triangle} $\phantom{.}$

In this case, $B$ is a small neighbourhood set of defect $2$,
and  $A$ and $C$ are small neighbourhood sets of defect $2$ or $3$. Lemma \ref{lm:estimate} shows
 that the number
of edges $r$ in $K$ satisfies $r \le 11$. Moreover, by the same lemma if $r=11$, then $|N_A'|= k-3 = |N_C'|$, and  if $r=10$ then $|N_A'| = k-3$ or $|N_C'| = k-3$. We will assume w.l.o.g. that $|N_A'| \le k-3$ whenever $r=10$. Moreover, if $r=11$ we will assume w.l.o.g. that there are at least as many edges from $B$ to $C$ as there are from $B$ to $A$.

\begin{lemma}\label{lm:322tri}
Each element of $x\in A \cup C$ is adjacent to at least 2 elements of $B$, and there is at least one edge from
$A$ to $C$.
\end{lemma}
\begin{pf}
The first part follows from property (\ref{NA}) of Lemma \ref{lm:ANA}. If there would be no edges between
$A$ and $C$, then $a_1,a_2$ could only be adjacent to the $3$ elements in $B$ and the $k-2$
elements in $N_A'$, leaving $|N_A|\le k+1$, for a contradiction.\qed
\end{pf}
Lemma \ref{lm:322tri} implies that $r \ge 9$, hence $K$ contains $9$, $10$, or $11$ edges.

\begin{lemma} There exists an element $x \in C$ that is adjacent to exactly one element of $A$.
\end{lemma}
\begin{pf} Lemma \ref{lm:322tri} together with the fact that $r$ satisfies $9\le r \le 11$ implies that there are $1$ to $3$ edges from $A$ to $C$. The statement of the Lemma is true unless there are exactly $2$ edges from $A$ to $C$ that share a vertex in $C$. Say w.l.o.g. that these are the edges $a_1-c_1-a_2$, so $c_2 \notin N_A$. Now, with the results of Lemma \ref{lm:322tri}, the $2$ edges between $A$ and $C$ require that $r\ge 10$, and hence $|N_A'|\le k-3$
by assumption. But then
$$|N_A| \le |N_A'|+|B|+|\{c_1\}| \le k+1,$$
contradicting $|N_A| \ge k+2$.
\qed
\end{pf}
Hence, we may assume that $c_1-a_1$, and that $c_1$ is non-adjacent to $a_2$. The following figure depicts the minimal amount of edges in $K$.

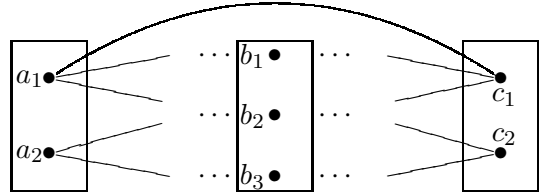
\[
\xy
(-5,5)*{}="x1";
(5,5)*{}="x4";
(-5,25)*{}="x2";
(5,25)*{}="x3";
"x1";"x2" **\dir{-};
"x2";"x3" **\dir{-};
"x3";"x4" **\dir{-};
"x4";"x1" **\dir{-};
(0,10)*{\bullet}="d2";
(0,20)*{\bullet}="d1";
(-2.6,10)*{a_2};
(-2.6,20)*{a_1};
(25,5)*{}="z1";
(35,5)*{}="z4";
(25,25)*{}="z2";
(35,25)*{}="z3";
"z1";"z2" **\dir{-};
"z2";"z3" **\dir{-};
"z3";"z4" **\dir{-};
"z4";"z1" **\dir{-};
(30,7)*{\bullet}="b2";
(30,23)*{\bullet}="b1";
(30,15)*{\bullet}="b3";
(27,7)*{b_3};
(27,23)*{b_1};
(27,15)*{b_2};
(38,15)*{\cdots};
(38,7)*{\cdots};
(38,23)*{\cdots};
(22,15)*{\cdots};
(22,7)*{\cdots};
(22,23)*{\cdots};
(46,14)*{}="cd3";
(46,17)*{}="cd31";
(45,7)*{}="cd2";
(45,23)*{}="cd1";
(15,14)*{}="xd3";
(15,17)*{}="xd31";
(16,7)*{}="xd2";
(16,23)*{}="xd1";
(55,5)*{}="w1";
(65,5)*{}="w4";
(55,25)*{}="w2";
(65,25)*{}="w3";
"w1";"w2" **\dir{-};
"w2";"w3" **\dir{-};
"w3";"w4" **\dir{-};
"w4";"w1" **\dir{-};
(60,10)*{\bullet}="c2";
(60,20)*{\bullet}="c1";
(60.5,12.3)*{c_2};
(60.5,17.5)*{c_1};
"c1";"cd1" **\dir{-};
"c1";"cd31" **\dir{-};
"c2";"cd2" **\dir{-};
"c2";"cd3" **\dir{-};
"d1";"xd1" **\dir{-};
"d1";"xd31" **\dir{-};
"d2";"xd2" **\dir{-};
"d2";"xd3" **\dir{-};
(0,0)*{}="jj";
(60,0)*{}="jjj";
"d1";"c1" **\crv{(30,40)} ?>* \dir{-};
(20,32)*{}="g1";
(70,32)*{}="g2";
(70,10)*{}="g3";
(20,-2)*{}="g4";
(67.5,-2)*{}="g5";
(67.5,20)*{}="g6";
\endxy
\]

By transitivity of $G$,
there exists a small neighbourhood set $D$ (the image of $B$ under some $g \in G$) of defect $2$ and size $3$ with $c_1\in D$. As $a_1 -c_1$, by Lemma \ref{lm:ANA}(c), there exists
$d \in D, d \ne c_1$ with $a_1-d$. Hence, $d \in  N_A' \cup B \cup \{c_2\}$. The following lemmas will examine these possibilities.

\begin{lemma} $d \notin B.$
\end{lemma}
\begin{pf} Assume otherwise, say that $d=b_1$. Consider the set $N=N(b_1) \cap N(c_1)$ with $|N|=k-2$. Then
$a_1$ is the only element in $N \cap A$, as $c_1 $ is not adjacent to $a_2$. The other elements of $N$ may not be in $B$ or $C$, as $b_1$ and $c_1$ are, and hence
are in singleton classes. 
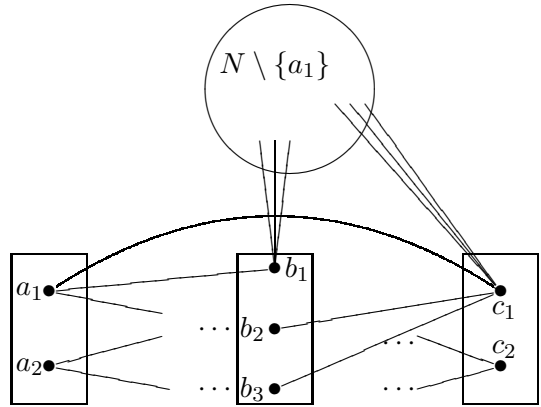
\[
\xy
(32,47)*\cir<35pt>{};
(15,40)*{}="q1";
(15,50)*{}="q2";
(45,40)*{}="q4";
(45,50)*{}="q3";
(30,50)*{N\setminus \{a_1\}};
(28,40)*{}="ghost1";
(30,40)*{}="ghost2";
(32,40)*{}="ghost3";
(40,45)*{}="ghost4";
(42,45)*{}="ghost5";
(38,45)*{}="ghost6";
(-5,5)*{}="x1";
(5,5)*{}="x4";
(-5,25)*{}="x2";
(5,25)*{}="x3";
"x1";"x2" **\dir{-};
"x2";"x3" **\dir{-};
"x3";"x4" **\dir{-};
"x4";"x1" **\dir{-};
(0,10)*{\bullet}="d2";
(0,20)*{\bullet}="d1";
(-2.6,10)*{a_2};
(-2.6,20)*{a_1};
(25,5)*{}="z1";
(35,5)*{}="z4";
(25,25)*{}="z2";
(35,25)*{}="z3";
"z1";"z2" **\dir{-};
"z2";"z3" **\dir{-};
"z3";"z4" **\dir{-};
"z4";"z1" **\dir{-};
(30,7)*{\bullet}="b2";
(30,23)*{\bullet}="b1";
(30,15)*{\bullet}="b3";
(27,7)*{b_3};
(33,23)*{b_1};
(27,15)*{b_2};
(46.5,13)*{\cdots};
(46.5,7)*{\cdots};
(22,15)*{\cdots};
(22,7)*{\cdots};
(49,14)*{}="cd3";
(46,17)*{}="cd31";
(49,7)*{}="cd2";
(45,23)*{}="cd1";
(15,14)*{}="xd3";
(15,17)*{}="xd31";
(16,7)*{}="xd2";
(16,23)*{}="xd1";
(55,5)*{}="w1";
(65,5)*{}="w4";
(55,25)*{}="w2";
(65,25)*{}="w3";
"w1";"w2" **\dir{-};
"w2";"w3" **\dir{-};
"w3";"w4" **\dir{-};
"w4";"w1" **\dir{-};
(60,10)*{\bullet}="c2";
(60,20)*{\bullet}="c1";
(60.5,12.3)*{c_2};
(60.5,17.5)*{c_1};
"c1";"b2" **\dir{-};
"c1";"b3" **\dir{-};
"c2";"cd2" **\dir{-};
"c2";"cd3" **\dir{-};
"d1";"xd31" **\dir{-};
"d1";"b1" **\dir{-};
"d2";"xd2" **\dir{-};
"d2";"xd3" **\dir{-};
"b1";"ghost1" **\dir{-};
"b1";"ghost2" **\dir{-};
"b1";"ghost3" **\dir{-};
"c1";"ghost4" **\dir{-};
"c1";"ghost5" **\dir{-};
"c1";"ghost6" **\dir{-};
(0,0)*{}="jj";
(60,0)*{}="jjj";
"d1";"c1" **\crv{(30,40)} ?>* \dir{-};
(20,32)*{}="g1";
(70,32)*{}="g2";
(70,10)*{}="g3";
(20,-2)*{}="g4";
(67.5,-2)*{}="g5";
(67.5,20)*{}="g6";
\endxy
\]

By Lemma \ref{lm:Nimage}, $b_1 f$ and $ c_1 f$ are non-adjacent. However, this is false, for a contradiction.
\qed
\end{pf}

\begin{lemma} $d \ne c_2.$
\end{lemma}
\begin{pf} Assume otherwise. Then $c_2-a_1$, and there are at least two edges between $A$ and $C$. Together with at least $4$ edges from $A$ to $B$, there are at most $5$ edges from $B$ to $C$. As $p_B=k+2$,
we have $B \subseteq N_C$, and with at most $5$ available edges, it follows that $N(c_1) \cap B \ne N(c_2) \cap B$.

Now let $e \notin \{c_1,c_2\} $ be the third element of $D$. We have that $B \subseteq N(c_1) \cup N(c_2)=N_D$. As $|N_D|=k+2$ and $|B|=3$, $e \in N_B$. 
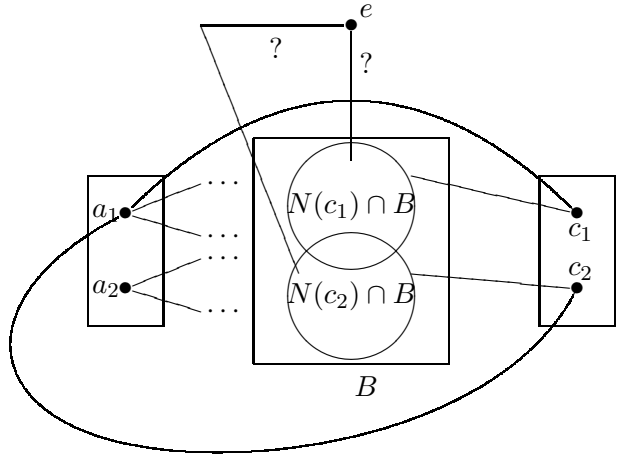
\[
\xy
(-5,5)*{}="x1";
(5,5)*{}="x4";
(-5,25)*{}="x2";
(5,25)*{}="x3";
"x1";"x2" **\dir{-};
"x2";"x3" **\dir{-};
"x3";"x4" **\dir{-};
"x4";"x1" **\dir{-};
(0,10)*{\bullet}="d2";
(0,20)*{\bullet}="d1";
(-2.6,10)*{a_2};
(-2.6,20)*{a_1};
(17,0)*{}="z1";
(43,0)*{}="z4";
(17,30)*{}="z2";
(43,30)*{}="z3";
"z1";"z2" **\dir{-};
"z2";"z3" **\dir{-};
"z3";"z4" **\dir{-};
"z4";"z1" **\dir{-};
(30,21)*\cir<25pt>{};
(30,21)*{N(c_1)\cap B};
(30,9)*\cir<25pt>{};
(30,9)*{N(c_2)\cap B};
(30,45)*{\bullet}="dot";
(32,47)*{e};
(30,27)*{}="dot1";
(23,12)*{}="dot2";
(10,45)*{}="dot3";
(38,12)*{}="dot4";
(38,25)*{}="dot5";
"dot";"dot1" **\dir{-};
"dot";"dot3" **\dir{-};
"dot2";"dot3" **\dir{-};
(32,40)*{?};
(20,42)*{?};
(32,-3)*{B};
(49,14)*{}="cd3";
(46,17)*{}="cd31";
(49,7)*{}="cd2";
(45,23)*{}="cd1";
(10,14)*{}="xd3";
(13,14)*{\cdots};
(10,17)*{}="xd31";
(13,17)*{\cdots};
(10,24)*{}="xd41";
(13,24)*{\cdots};
(10,7)*{}="xd2";
(13,7)*{\cdots};
(10,23)*{}="xd1";
(55,5)*{}="w1";
(65,5)*{}="w4";
(55,25)*{}="w2";
(65,25)*{}="w3";
"w1";"w2" **\dir{-};
"w2";"w3" **\dir{-};
"w3";"w4" **\dir{-};
"w4";"w1" **\dir{-};
(60,10)*{\bullet}="c2";
(60,20)*{\bullet}="c1";
(60.5,12.3)*{c_2};
(60.5,17.5)*{c_1};
"d1";"xd31" **\dir{-};
"d1";"xd41" **\dir{-};
"d2";"xd2" **\dir{-};
"d2";"xd3" **\dir{-};
(0,0)*{}="jj";
(60,0)*{}="jjj";
"d1";"c2" **\crv{(-20,10) & (-20,-20) & (50,-10)} ?>* \dir{-};
"d1";"c1" **\crv{(30,50)} ?>* \dir{-};
(20,32)*{}="g1";
(70,32)*{}="g2";
(70,10)*{}="g3";
(20,-2)*{}="g4";
(67.5,-2)*{}="g5";
(67.5,20)*{}="g6";
"c2";"dot4" **\dir{-};
"c1";"dot5" **\dir{-};
\endxy
\]

Now, $N(e) \cap B$ must
differ from one of $N(c_1) \cap B$, $N(c_2) \cap B$. This contradicts Lemma \ref{lm:defect2b}(\ref{mapclique}), for $ef \ne c_1 f= c_2 f$.
\qed
\end{pf}
\begin{lemma} $d \notin N_A'.$
\end{lemma}
\begin{pf} Assume otherwise. By Lemma \ref{lm:defect2b}(\ref{mapclique}), $N(d) \cap A=N(c_1) \cap A=\{a_1\}$, and so $a_2 \notin N(d)$. This implies that $a_2$ must have at least $k -\left(|N_A'|-1\right)$ neighbours in $K$.

Now, if $r=9$, then $|N_A'|\le k-2$, and so $a_2$ requires at least $3$ neighbors in $K$. However, Lemma \ref{lm:322tri} accounts for all $9$ edges in $K$, showing that $a_2$ has exactly $2$ neighbors in $K$ (recall that the edge from $A$ to $C$ was assumed to be $a_1-c_1$). This excludes
the case $r=9$.

If $r\ge 10$, then  $|N_A'|= k-3$, and so $a_2$ requires at least $4$ neighbors in $K$, which must be the elements of $B \cup \{c_2\}$. With $3$ edges from $a_2$ to $B$, $a_2-c_2$, $a_1-c_1$, $2$ edges from $a_1$ to $B$, and
$4$ edges between $B$ and $C$, we see that $r=11$.
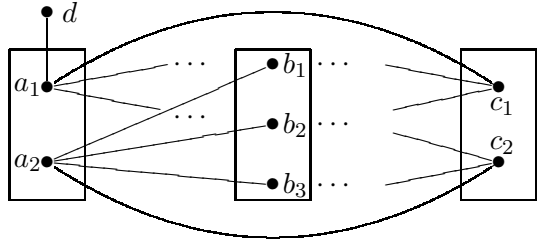
\[
\xy
(-5,5)*{}="x1";
(5,5)*{}="x4";
(-5,25)*{}="x2";
(5,25)*{}="x3";
"x1";"x2" **\dir{-};
"x2";"x3" **\dir{-};
"x3";"x4" **\dir{-};
"x4";"x1" **\dir{-};
(0,10)*{\bullet}="d2";
(0,20)*{\bullet}="d1";
(0,30)*{\bullet}="d";
(3,30)*{d};
"d";"d1" **\dir{-};
(-2.6,10)*{a_2};
(-2.6,20)*{a_1};
(25,5)*{}="z1";
(35,5)*{}="z4";
(25,25)*{}="z2";
(35,25)*{}="z3";
"z1";"z2" **\dir{-};
"z2";"z3" **\dir{-};
"z3";"z4" **\dir{-};
"z4";"z1" **\dir{-};
(30,7)*{\bullet}="b2";
(30,23)*{\bullet}="b1";
(30,15)*{\bullet}="b3";
(33,7)*{b_3};
(33,23)*{b_1};
(33,15)*{b_2};
(38,15)*{\cdots};
(38,7)*{\cdots};
(38,23)*{\cdots};
(19,16)*{\cdots};
(19,23)*{\cdots};
(46,14)*{}="cd3";
(46,17)*{}="cd31";
(45,7)*{}="cd2";
(45,23)*{}="cd1";
(15,14)*{}="xd3";
(15,17)*{}="xd31";
(16,7)*{}="xd2";
(16,23)*{}="xd1";
(55,5)*{}="w1";
(65,5)*{}="w4";
(55,25)*{}="w2";
(65,25)*{}="w3";
"w1";"w2" **\dir{-};
"w2";"w3" **\dir{-};
"w3";"w4" **\dir{-};
"w4";"w1" **\dir{-};
(60,10)*{\bullet}="c2";
(60,20)*{\bullet}="c1";
(60.5,12.3)*{c_2};
(60.5,17.5)*{c_1};
"c1";"cd1" **\dir{-};
"c1";"cd31" **\dir{-};
"c2";"cd2" **\dir{-};
"c2";"cd3" **\dir{-};
"d1";"xd1" **\dir{-};
"d1";"xd31" **\dir{-};
"d2";"b1" **\dir{-};
"d2";"b2" **\dir{-};
"d2";"b3" **\dir{-};
(0,0)*{}="jj";
(60,0)*{}="jjj";
"d1";"c1" **\crv{(30,40)} ?>* \dir{-};
"d2";"c2" **\crv{(30,-10)} ?>* \dir{-};
(20,32)*{}="g1";
(70,32)*{}="g2";
(70,10)*{}="g3";
(20,-2)*{}="g4";
(67.5,-2)*{}="g5";
(67.5,20)*{}="g6";
\endxy
\]

However, for the case that $r=11$, we assumed that there are at least as many edges from $B$ to $C$ as there are from $B$ to $A$. Our final configuration violates this assumption, for a contradiction.

\qed
\end{pf}
We have excluded every possible location for $d$. Therefore, Theorem \ref{th:n-4} holds for transformations $f$ of kernel type $(3,2,2,$ $1,\dots,1)$.

\subsection{Maps of kernel type $(2,2,2,2,1, \dots,1)$}

Let $A,B,C,D$ be the non-singleton kernel classes of $f$, and let $A=\{a_1,a_2\}$,  $B=\{b_1,b_2\}$, $C=\{c_1,c_2\}$, $D=\{d_1,d_2\}$.

For each kernel class $Z$ with image $z'$, $p_Z\ge k+2$ implies that $z'$ must be adjacent to
at least $2$  other images of non-singleton kernel classes. Hence the induced subgraph on $Kf$ must have $6$, $5$, or $4$ edges, and in the last case, these must form a $4$-cycle.

Throughout, $g$ will denote the transformation $(a_1\, a_2)(b_1\, b_2)(c_1\, c_2)(d_1\, d_2)$. As noted in Subsection~\ref{smallmindegree}, we are done if
we can establish that $g$ is an automorphism of $\Gam'$.

\paragraph{The image of $K$ has $4$ edges arranged in a cycle} $\phantom{.}$

We may suppose that the images of the non-singleton kernel classes are  $a_1 f-b_1 f- c_1 f- d_1 f-a_1 f$.
In this case each non-trivial kernel class $Z$ satisfies $p_Z=k+2$, and is hence a small neighbourhood set of
defect $2$. Hence $Z_1 \subseteq N_{Z_2}$ for every pair $(Z_1,Z_2)$ of adjacent kernel classes. In particular, there are at least two edges between each such pair.

Let $r$ be the number of edges in $K$. By Lemma \ref{lm:estimatelow} and Lemma \ref{lm:estimate}, we have $8 \le r \le 10$.

\subparagraph{$K$ contains $8$ edges} $\phantom{.}$

Here there are exactly two edges between each pair $(Z_1,Z_2)$ of adjacent kernel classes.
Now $Z_1 \subseteq N_{Z_2}$ and $Z_2 \subseteq N_{Z_1}$ is only possible if the two edges between $Z_1$ and $Z_2$ have
disjoint vertices. The only two possible configurations are depicted in Figures \ref{fi:im4ed8a} and \ref{fi:im4ed8b}.
\begin{figure}[h]
\[
\xy
(-25,5)*{}="0b";
(-15,5)*{}="0c";
(-25,25)*{}="0f";
(-15,25)*{}="0g";
"0b";"0f" **\dir{-};
"0b";"0c" **\dir{-};
"0c";"0g" **\dir{-};
"0g";"0f" **\dir{-};
(-20,20)*{\bullet}="0e";
(-20,10)*{\bullet}="0d";
(-5,5)*{}="b";
(5,5)*{}="c";
(-5,25)*{}="f";
(5,25)*{}="g";
"b";"f" **\dir{-};
"b";"c" **\dir{-};
"c";"g" **\dir{-};
"g";"f" **\dir{-};
(15,5)*{}="bS";
(25,5)*{}="cS";
(15,25)*{}="fS";
(25,25)*{}="gS";
"bS";"fS" **\dir{-};
"bS";"cS" **\dir{-};
"cS";"gS" **\dir{-};
"gS";"fS" **\dir{-};
(35,5)*{}="bT";
(45,5)*{}="cT";
(35,25)*{}="fT";
(45,25)*{}="gT";
"bT";"fT" **\dir{-};
"bT";"cT" **\dir{-};
"cT";"gT" **\dir{-};
"gT";"fT" **\dir{-};
(0,20)*{\bullet}="e";
(0,10)*{\bullet}="d";
(-20,0)*{}="a";
(40,10)*{\bullet}="d2";
(40,0)*{}="a2";
(20,10)*{\bullet}="d1";
(20,20)*{\bullet}="a1T";
(40,20)*{\bullet}="a2T";
(-20,30)*{}="dT";
(40,30)*{}="d2T";
"d";"0d" **\dir{-};
"a";"0d" **\dir{-};
"a2";"a" **\dir{-};
"d2";"d1" **\dir{-};
"d1";"d" **\dir{-};
"0e";"dT" **\dir{-};
"0e";"e" **\dir{-};
"dT";"d2T" **\dir{-};
"e";"a1T" **\dir{-};
"a1T";"a2T" **\dir{-};
"d2T";"a2T" **\dir{-};
"a2";"d2" **\dir{-};
\endxy
\]
\caption{One of the two configuration with $8$ edges}\label{fi:im4ed8a}
\end{figure}
\begin{figure}[h]
\[
\xy
(-25,5)*{}="0b";
(-15,5)*{}="0c";
(-25,25)*{}="0f";
(-15,25)*{}="0g";
"0b";"0f" **\dir{-};
"0b";"0c" **\dir{-};
"0c";"0g" **\dir{-};
"0g";"0f" **\dir{-};
(-20,20)*{\bullet}="0e";
(-20,10)*{\bullet}="0d";
(-5,5)*{}="b";
(5,5)*{}="c";
(-5,25)*{}="f";
(5,25)*{}="g";
"b";"f" **\dir{-};
"b";"c" **\dir{-};
"c";"g" **\dir{-};
"g";"f" **\dir{-};
(15,5)*{}="bS";
(25,5)*{}="cS";
(15,25)*{}="fS";
(25,25)*{}="gS";
"bS";"fS" **\dir{-};
"bS";"cS" **\dir{-};
"cS";"gS" **\dir{-};
"gS";"fS" **\dir{-};
(35,5)*{}="bT";
(45,5)*{}="cT";
(35,25)*{}="fT";
(45,25)*{}="gT";
"bT";"fT" **\dir{-};
"bT";"cT" **\dir{-};
"cT";"gT" **\dir{-};
"gT";"fT" **\dir{-};
(0,20)*{\bullet}="e";
(0,10)*{\bullet}="d";
(-20,0)*{}="a";
(40,10)*{\bullet}="d2";
(40,0)*{}="a2";
(20,10)*{\bullet}="d1";
(20,20)*{\bullet}="a1T";
(40,20)*{\bullet}="a2T";
(-20,30)*{}="dT";
(40,30)*{}="d2T";
"d";"0d" **\dir{-};
"a";"0d" **\dir{-};
"a2";"a" **\dir{-};
"d2";"d1" **\dir{-};
"d1";"e" **\dir{-};
"0e";"dT" **\dir{-};
"0e";"e" **\dir{-};
"dT";"d2T" **\dir{-};
"d";"a1T" **\dir{-};
"a1T";"a2T" **\dir{-};
"d2T";"a2T" **\dir{-};
"a2";"d2" **\dir{-};
\endxy
\]
\caption{One of the two configuration with $8$ edges}\label{fi:im4ed8b}
\end{figure}

It is now easy to check that $g$ is an automorphism of $\Gam'$, for a contradiction.
\subparagraph{$K$ contains $9$ edges} $\phantom{.}$

We may assume that $A$ and $B$ are the unique non-singleton kernel classes that have $3$ edges between
them, and that
$b_1-a_1-b_2-a_2$. However in this case, we have $N_B' \subseteq N(b_1)$, which implies the CME$(a_1-b_2,b_1)$, for a contradiction.
\subparagraph{$K$ contains $10$ edges} $\phantom{.}$

Suppose first that we have two kernel classes that have only three edges between them, say $A$ and $B$ with edges $b_1-a_1-b_2-a_2$.
By the number of available edges, at least one of $b_1$, $a_2$
is a vertex of only two edges from within $K$. Hence either $N_B' \subseteq N(b_1)$
or  $N_A' \subseteq N(a_2)$, and so we have the CME$(a_1-b_2,b_1)$ or CME$(a_1-b_2,a_2)$, as in the case that $K$ contains $9$ edges.

\begin{figure}[h]
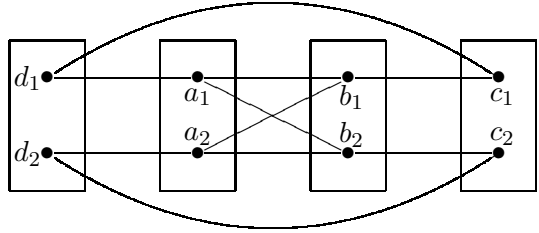

\[
\xy
(-5,5)*{}="x1";
(5,5)*{}="x4";
(-5,25)*{}="x2";
(5,25)*{}="x3";
"x1";"x2" **\dir{-};
"x2";"x3" **\dir{-};
"x3";"x4" **\dir{-};
"x4";"x1" **\dir{-};
(0,10)*{\bullet}="d2";
(0,20)*{\bullet}="d1";
(-2.6,10)*{d_2};
(-2.6,20)*{d_1};
(15,5)*{}="y1";
(25,5)*{}="y4";
(15,25)*{}="y2";
(25,25)*{}="y3";
"y1";"y2" **\dir{-};
"y2";"y3" **\dir{-};
"y3";"y4" **\dir{-};
"y4";"y1" **\dir{-};
(20,10)*{\bullet}="a2";
(20,20)*{\bullet}="a1";
(20,12.3)*{a_2};
(20,17.5)*{a_1};
(35,5)*{}="z1";
(45,5)*{}="z4";
(35,25)*{}="z2";
(45,25)*{}="z3";
"z1";"z2" **\dir{-};
"z2";"z3" **\dir{-};
"z3";"z4" **\dir{-};
"z4";"z1" **\dir{-};
(40,10)*{\bullet}="b2";
(40,20)*{\bullet}="b1";
(40.5,12.5)*{b_2};
(40.5,17.5)*{b_1};
(55,5)*{}="w1";
(65,5)*{}="w4";
(55,25)*{}="w2";
(65,25)*{}="w3";
"w1";"w2" **\dir{-};
"w2";"w3" **\dir{-};
"w3";"w4" **\dir{-};
"w4";"w1" **\dir{-};
(60,10)*{\bullet}="c2";
(60,20)*{\bullet}="c1";
(60.5,12.3)*{c_2};
(60.5,17.5)*{c_1};
"d2";"a2" **\dir{-};
"a1";"b1" **\dir{-};
"a2";"b2" **\dir{-};
"d1";"a1" **\dir{-};
"a1";"b2" **\dir{-};
"a2";"b1" **\dir{-};
"b1";"c1" **\dir{-};
"b2";"c2" **\dir{-};
"d2";"c2" **\crv{(30,-10)} ?>* \dir{-};
"d1";"c1" **\crv{(30,40)} ?>* \dir{-};
(20,32)*{}="g1";
(70,32)*{}="g2";
(70,10)*{}="g3";
(20,-2)*{}="g4";
(67.5,-2)*{}="g5";
(67.5,20)*{}="g6";
\endxy
\]
\caption{The remaining configuration for $K$ with $10$ edges}\label{fi:im4ed10}
\end{figure}

So assume instead that there are two kernel classes with $4$ edges between them, say $A$ and $B$; then we may assume that all edges within $K$ are  $a_1-b_1-a_2-b_2-a_1$, $b_1-c_1-d_1-a_1$, and $b_2-c_2-d_2-a_2$ (see Figure \ref{fi:im4ed10}).
We have $N_C' \subseteq N(c_1) \cap N(c_2)$ and $N_D' \subseteq  N(d_1) \cap N(d_2)$. Moreover $N_B'$ consists of $k-4$ elements that are
 adjacent to both $b_1$ and $b_2$, one element $b_1'$ adjacent to $b_1$ but not $b_2$, and one element $b_2'$ adjacent to $b_2$ but not $b_1$.

\begin{lemma} $|N(b_1') \cap N(c_1)|=k-2$.\label{lm:b1pc1}
\end{lemma}
\begin{pf}
By Lemma \ref{lm:smnd2s2} applied to $y=c_1, y'=b_1, \bar y=d_1$, there exist $z \in \Gam'$ with $|N(c_1) \cap N(z)|=k-2$, such that $z$ is
adjacent to one of  $b_1,d_1$. We want to narrow the location of $z$.

As $N(z)\cap \{b_1,d_1\} \ne \emptyset$,
$z \in \{a_1,a_2\} \cup (N_B' \setminus\{b_2'\}) \cup N_D'$.
If $z \in \{a_1,a_2\}$ then $N(z) \cap B = B \ne  \{b_1\} =N(c_1 ) \cap B$, contradicting Lemma \ref{lm:defect2b}. Similarly, if $z \in N_D'$ then $N(z) \cap D = D\ne \{d_1\}= N(c_1) \cap D$,
 and if $z \in N_B' \setminus \{b_1',b_2'\}$ then $N(z) \cap B = B \ne  \{b_1\} =N(c_1 ) \cap B$. Hence $z=b_1'$.\qed
\end{pf}
\begin{lemma} \label{lm:nhdb1'} $N(b_1')=\{a_1,a_2,b_1\} \cup H$ where $H \subseteq N_C'$.
\end{lemma}

\begin{pf} $b_1'$ must be in every $r$-clique containing $a_1-b_1$,
for otherwise we obtain a CME$(a_1-b_1,b_2)$. Hence $a_1 \in N(b_1')$. Similarly, $a_2 \in N(b_1')$
to avoid a CME($a_2-b_1,b_2)$.

Hence $N(b_1') \setminus N(c_1)=\{a_1,a_2\}$. By Lemma \ref{lm:b1pc1},
 all remaining neighbours of $b_1'$ are in $N(c_1)$. One of those elements is $b_1$. If
$d_1 \in N(b_1')$, then
$N(b_1') \cap A = A \ne  \{a_1\} =N(d_1) \cap A$, contradicting Lemma \ref{lm:defect2b}. Hence
$N(b_1) \setminus \{a_1,a_2,b_1\} \subseteq  N_C'$. \phantom{PH}\hfill\qed
\end{pf}
\begin{lemma} There exists $x \in \Gam'$ such that $|N(b_1) \cap N(x)|=k-2$, $x$ is adjacent to $b_1'$, and $x$ is not adjacent to $c_1$.
\end{lemma}
\begin{pf}
By Lemma \ref{lm:smnd2s2} applied to $y=b_1, y'=b_1', \bar y=c_1$, there exist $x \in \Gam'$ with $|N(b_1) \cap N(x)|=k-2$, such that $x$ is
adjacent to one of  $b_1',c_1$.

If $x$ is adjacent to $c_1$ then either $x=d_1$ or $x \in N_C'$. Now if $x=d_1$ then $N(b_1) \cap A = A \ne  \{a_1\} =N(d_1) \cap A$, contradicting Lemma \ref{lm:defect2b}. Similarly,  if $x \in N_C'$, then
 $N(x) \cap C = C \ne  \{c_1\} =N(b_1 ) \cap C$. Hence $x \in N(b_1') \setminus N(c_1)$.\qed
\end{pf}
By Lemma \ref{lm:nhdb1'}, we get that $x \in A$. This implies that $b_1-x$,  contradicting Lemma \ref{lm:defect2b}(\ref{nonadj}).

\paragraph{The image of $K$ has $5$ edges} $\phantom{.}$

We may assume that the edges in the image of $K$ are  $a_1 f- b_1 f - c_1 f- d_1 f - a_1 f-c_1f$. Hence $B$ and $D$ are small neighbourhood classes of defect $2$ and $A$ and $C$ are small neighbourhood classes of defect $2$ or $3$. As $p_B=p_D=k-2$, there are at least $2$ edges between each kernel class pair in $\{A,C\} \times \{B,D\}$.

Let $r$ be the number of edges in $K$; by Lemma \ref{lm:estimatelow} and Lemma \ref{lm:estimate}, we
obtain $10\le r \le 13$.
Moreover,
$r=13 $ implies that $|N_A'|= k-4 = |N_C'|$, and $r=12$ implies that $|N_A'| = k-4$ or $|N_C'| = k-4$.

\begin{lemma}\label{lm:detour}
Let $z \in B \cup D$. Suppose that $|N(z) \cap K|=2$. Then
$|N(z) \cap A|=1=|N(z)\cap C|$.
\end{lemma}
\begin{pf}
Suppose otherwise, say w.l.o.g that $N(b_1) \cap K=A$.

We  claim that  one element $x \in A$ satisfies $N_A' \subseteq N(x)$.
If $r \le 11$, then at most $7$ edges have a vertex in $A$, as at least $4$ edges lie between $C,B$ and $C,D$. Thus one of $a_1,a_2$ must have $k-3=|N_A'|$ neighbours outside of $K$.

If $r \ge 12$, then $|N_A'|= k-4$ or
$|N_C'| = k-4$. However $b_1 \notin N_C$, and so $N_C \cap K$ has at most $5$ elements. As $|N_C| \ge k+2$, it follows that $|N_C'| = k-3$, and so $|N_A'| = k-4$.
Because $r \le 13$  at most $9$ edges have a vertex in $A$, so one of $a_1,a_2$ must have $k-4$ neighbours outside of $K$. 

In either case $N_A' \subseteq N(x)$ for some $x \in A$, say for $a_1$.
However, we now have a CME($b_1-a_2,a_1)$,  for a contradiction. So $|N(b_1) \cap A|=1$, and thus $|N(b_1)\cap C|=1$.
\qed
\end{pf}
\begin{lemma}$\Gam'$ has at least 10 edges that lie between the pairs of kernel classes from $\{A,C\} \times \{B,D\}$. \label{lm:5imageGen}
\end{lemma}
\begin{pf}
Assume to the contrary that there are at most $9$ edges between the pairs of kernel classes from $\{A,C\} \times \{B,D\}$. We will construct a contradiction to Lemma \ref{lm:Nimage}.

As there are at least two edges between the pairs in $\{A,C\} \times \{B,D\}$, each pair has either $2$ or $3$ edges between them,
with at most one case of $3$ edges. We may assume that the exceptional pair in the case of $3$ edges  is $(C,D)$. Applying Lemma \ref{lm:detour} to the $3$ or $4$
vertices $z \in B \cup D$ that have exactly two neighbours in $K$, we see that if there are two edges between any pair $(Y,Z)$ of kernel classes, those edges have disjoint vertices.

Hence, w.l.o.g. we may assume that we have the edges $b_1-a_1-d_1$ and $b_2-a_2-d_2$. In case that there are $3$ edges between $C$ and $D$, we may further assume that $d_1$ is the unique vertex in $D$ with $3$ neighbours in $K$.
Applying Lemma \ref{lm:smnd2s2} with $y=a_1$, $y'=b_1$, $\bar y=d_1$, we see that there is a $z$ such that $|N(a_1)\cap N(z)|=k-2$, with $z$ adjacent to $b_1$ or $d_1$.

We claim that $z \in C$.
As $z \in N(b_1)\cup N(d_1)$, we have $z \in C \cup N_B' \cup N_D'$. Now for all $w \in N_D'$, $w \in N(d_2)$ as $d_2$ has only two neighbours in $K$. Hence
$N(w) \cap D \ne \{d_1\} = N(a_1) \cap D$, and so $z \notin N_D'$ by Lemma \ref{lm:defect2b}(\ref{mapclique}). An analogous argument show that $z \notin N_B'$, and so $z \in C$.

Let $N= N(a_1) \cap N(z)$. We claim that for every non-singleton kernel class $Z$ of $f$, $|N \cap Z| \le 1$.
$N \cap B \subseteq N(a_1) \cap B=\{b_1\}$ and   $N \cap D \subseteq N(a_1) \cap D=\{d_1\}$, so the claim holds for $Z=B$ and $Z=D$. Moreover, $N$ does not have any elements in $A$ or $C$, as $a_1 \in A, z \in C$.

Hence Lemma \ref{lm:Nimage} is applicable to $N$. By the lemma $a_1 f$ and $z f$ are non-adjacent. However,
we have that $a_1 f-c_1 f=z f$, as $z \in C$ , for a contradiction.
$\phantom{YY}$\hfill \qed
\end{pf}

\subparagraph{$K$ contains $10$ or $11$ edges} $\phantom{.}$

By Lemma \ref{lm:5imageGen}, in these cases there is at most one edge between $A$ and $C$. Our next Lemma shows that this is not possible, for a contradiction.
\begin{lemma}\label{lm:k11diag2} If $r\le 11$, there are at least two edges from $A$ to $C$.
\end{lemma}
\begin{pf}
At least one edge must cross from $A$ to $C$, for otherwise not all elements in
$A \cup C$ could have $3$ neighbours in $K$.

Assume that there is only one edge between $A$ and $C$.
 As at least $6$ edges go from $A$ to $K\setminus A$, there must be at least $5$ from $A$ to $B \cup D$, and
by symmetry at least $5$ edges from $C$ to $B \cup D$. This accounts for the maximum $11$ edges. Hence there are exactly $5$ edges from $A$ to $B \cup D$.

W.l.o.g. we may assume that there are $3$ edges from $A$ to $D$, say $a_1-d_1-a_2-d_2$, and $2$ edge from $A$ to $B$. The two edges from $A$ to $B$ must be adjacent to different elements of $A$ as $A \subseteq N_B$.
 This implies that the edge between $A$ and $C$ is adjacent to
$a_1$, and hence $N(a_2) \cap C= \emptyset$. Moreover, $N(a_1) \cap \bar K=N_A'$,
as $a_1$ has only three neighbours in $K$.

However, we now obtain CME$(a_2-d_1,a_1)$ for a contradiction. Hence there are at least two edges
 between $A$ and $C$.
\qed \end{pf}

\subparagraph{$K$ contains $12$ edges} $\phantom{.}$

In this case $|N_A'|=k-4$ or $|N_C'|=k-4$, say $|N_A'|=k-4$. Hence at least $8$  edges go from $A$ to $ K \setminus A$, while at least $6$ edges go from $C$ to $K \setminus C$. With $r=12$ this implies that at least
$2$ edges lie between $A$ and $C$.
 With Lemma \ref{lm:5imageGen}, we see that there
 are exactly $2$ edges between $A$ and $C$. As $C$ needs to be contained in $N_B$ and $N_D$, there exactly $2$ edges each between $(C, B)$ and $(C,D)$. This leaves $6$ edges between $(A,B)$ and $(A,D)$, and
  all edges are accounted for. Hence $a_1,a_2$ are both adjacent to exactly $4$ elements in $K$, and thus $N_A' \subseteq N(a_1)\cap N(a_2)$.

Assume first that there are $3$ edges between each of these pairs, where we may assume that $b_1-a_1-b_2-a_2$.
We have CME$(a_1-b_2,a_2)$, unless there is an element in $C$ (which we may assume to be $c_1$)
such that $b_2-c_1 -a_1$, and that $a_2$ is not adjacent to $c_1$. This implies that the second edge between $B$ and $C$ is $b_1-c_2$, and so in particular $c_2 \notin N(b_2)$.
But then $N(a_2) \cap N(b_2) \cap C= \emptyset$, and we obtain CME$(a_2-b_2, a_1)$,
for a contradiction.

Up to symmetry, the only remaining option is that there are $4$ edges between
 $A$ and $B$, and $2$ edges between $(A,D)$.
We obtain a CME$(a_1-b_1, a_2)$,  unless one element of $C$, say $c_1$, satisfies $a_1-c_1-b_1$ and $c_1 \notin N(a_2)$.
However, we now obtain a CME$(a_1-b_2, a_2)$, unless there exists $x \in C$ satisfying $a_1-x-b_2$ and that $x \notin N(a_2)$. $x \ne c_1$, for otherwise $c_2 \notin N_B$, as there are only two edges from $C$ to $B$.
 Hence $x=c_2$ and $N(a_2) \cap C= \emptyset$.
Finally, we obtain the CME$(a_2-b_1, a_1)$, for a contradiction.

Hence we can exclude the possibility that $K$ has $12$ edges.

\subparagraph{$K$ contains $13$ edges} $\phantom{.}$

By Lemma \ref{lm:estimate}, $|N_A'|=|N_C'|=k-4$. Hence each element of $A \cup C$ has at least $4$ neighbours in $K$, and as $r=13$ this is only possible if there are at least $3$ edges from $A$ to $C$. In fact, Lemma \ref{lm:5imageGen} show that there are exactly $3$ edges between $A$ and $C$, which in turn implies that each $x \in A \cup C$ has exactly $4$ neighbours in $K$.
 This implies that $N_A' \subseteq N(a_1) \cap N(a_2)$.

Up to symmetry, we may assume that there are $3$ edges from $A$ to $B$, and $2$ edges from $A$ to $D$, say that $b_1-a_1-b_2-a_2$. As there are 3 edges between $A$ and $C$ one of $a_1,a_2$ is adjacent to both elements in $C$. This must be $a_2$, for otherwise $a_1$ has $4$ neighbours in $B \cup C$ and could not be in $N_D$. So $c_1-a_2-c_2$.
But then we have a CME$(a_1-b_2, a_2)$ for a final contradiction.

\paragraph{The image of $K$ has $6$ edges} $\phantom{.}$

Now let $r$ be the number of edges between the elements of $K$. By Lemmas \ref{lm:estimatelow} and
\ref{lm:estimate} we get $12\le r \le 16$.
Moreover, by Lemma  \ref{lm:estimate}, if $p=r-12$, there are at least $p$ non-singleton kernel classes $Z$ for which $|N_Z'|\le k-4$.

Conversely, if there  are $p$ non-singleton kernel classes $Z$ for which $|N_Z'|\le k-4$, there are at least $8$ edges from each such $Z$ to $K \setminus Z$ and at least $6$ edges from any other class $Y$ to $K \setminus Y$.
This requires at least $(8p+6(4-p))/2=12+p=r$ edges. Hence if there are $12+p$ edges, there are exactly $p$  kernel classes $X$ for which $|N_X'|=k-4$, and exactly $4-p$ kernel classes with
$|N_Z'|=k-3$. As this accounts for all edges, we have proved the
following lemma.

\begin{lemma}\label{lm:fulledges}
Let $Z$ be a non-singleton kernel class, and $x \in Z$. If $|N_Z'|=k-3$, then $x$ has exactly $3$ neighbours in $K$. If $|N_Z'|=k-4$, then $x$ has exactly $4$ neighbours in $K$. In particular, $N_Z' \subseteq N(x)$.
\end{lemma}

\begin{lemma}\label{lm:k-3class}
Let $x \in Z$, where $Z$ is a kernel class with $|N_Z'|=k-3$. Then all three neighbours of $x$ in $K$ lie in different kernel classes.
\end{lemma}
\begin{pf}
Suppose otherwise, say w.l.o.g. that $x=b_1$, and that $a_1-b_1-a_2$.  Then we have CME$(a_1-b_1, a_2)$,  unless there exists $x \in C \cup D$ satisfying $a_1-x-b_1$ and that $x \notin N(a_2)$.
 Hence $a_1,a_2, x$ account for all neighbours of $b_1$ in $K$. But now we have CME$(a_2-b_1,a_1)$, for a contradiction.
 \qed
\end{pf}

Note that if $|N_Z'|=k-4$, then $|N_Z| = k+2$, and so $Z$ is a small neighbourhood set of defect $2$.

\subparagraph{$K$ contains $12$ edges} $\phantom{.}$

 Then $|N_X'|=k-3$ for all $X$ and by Lemma \ref{lm:k-3class}, every element of $K$ has exactly $3$ neighbours in $K$, all from different kernel classes. This implies that there are exactly
 $2$ edges between each pair of kernel classes, and that these edges have disjoint vertices.
 It follows that $g$ is an automorphism, and the result follows.

\subparagraph{$K$ contains $13, 14, $ or $15$ edges} $\phantom{.}$

In this case, we have kernel classes $Y, Z$ such that $|N_Y'|=k-4$, $|N_Z'|=k-3$. Note that $Y$ is a small neighbourhood set of defect $2$.

Assume w.l.o.g. that $Z=A$, then by Lemma \ref{lm:k-3class}, we may
assume that $N(a_1)\cap K=\{b_1,c_1,d_1\}$. $|N(a_1) \cap N(a_2)\cap K|\le 1$, for otherwise $|N_A| < k+2$. Thus we
may further assume that $b_2-a_2-c_2$. By Lemma  \ref{lm:k-3class}, $a_2$ has no additional
neighbours in $B \cup C$.

Now applying Lemma \ref{lm:smnd2s2} with $y=a_1, y'=b_1, \bar y=c_1$ there exists $z \in N(b_1) \cup N(c_1)$ with $|N(z) \cap N(a_1)|=k-2$. $z \ne a_2$, as $a_2$ is not adjacent to $b_1$ or $c_1$.
If $w \in N_B'$, then $N(w)\cap B=B \ne \{b_1\}=N(a_1) \cap B$, and so $w\ne z$ by Lemma \ref{lm:defect2b}(\ref{mapclique}). Analog, we get that $z \notin N_C'$.
It follows that $z \in B \cup C \cup D$.

Now consider $N= N(a_1) \cap N(z)$. As $N(a_1)$ intersects every kernel class in at most one point, the same holds for $N$. By Lemma \ref{lm:Nimage}, $a_1 f $ and $z f$ are non-adjacent.
However, as $z \in B \cup C \cup D$, this is false, for a contradiction.

\subparagraph{$K$ contains $16$ edges} $\phantom{.}$

Here  $|N_Z'|=k-4$ for all non-singleton kernel classes $Z$. As $|N_Z| \ge k+2$, this implies that $K \setminus Z \subseteq N_Z$. It follows that if there are exactly two edges between a pair of kernel classes,
those edges have disjoint vertices.
By Lemma \ref{lm:fulledges}, each element of $K$ has exactly four neighbours in $K$. Up to symmetry, there are two possibilities:
\begin{enumerate}
  \item There are $4$ edges between $A$ and $B$, $4$ edges between $C$ and $D$, and $2$ edges each between the other pairs of kernel classes;
  \item There are $2$ edges between $A$ and $B$, $2$ edges between $C$ and $D$, and $3$ edges each between the other pairs of kernel classes;\label{l:finalcase}
\end{enumerate}
In the first case, it is easy to see that $g$ is a graph automorphism, as the edges between pairs of classes other than $(A,B)$ and $(C,D)$ have disjoint vertices. So assume we are in the situation  (\ref{l:finalcase}).

 We may assume that the three edges between $A$ and $C$ are $c_1-a_1-c_2-a_2$. Hence $c_2$ has two neighbours in $A$, one neighbour in $D$, and thus one neighbour in $B$, which we may assume to be $b_2$. Similarly,
 $c_1$ has two neighbours in $B$, and we get the edges $b_1-c_1-b_2-c_2$ between $B$ and $C$. Continuing in this fashion, we get the edges   $d_1-b_1-d_2-b_2$ and  $a_1-d_1-a_2-d_2$.

\begin{figure}[h]
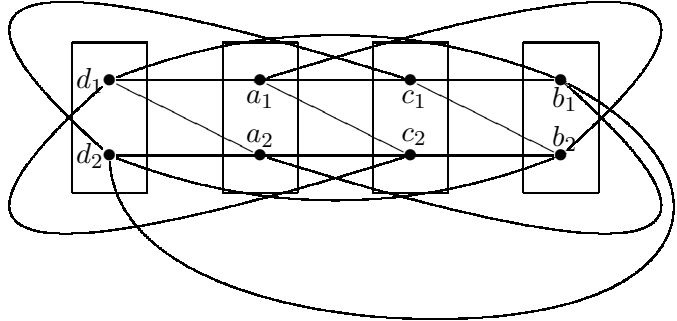

\[
\xy
(-5,5)*{}="x1";
(5,5)*{}="x4";
(-5,25)*{}="x2";
(5,25)*{}="x3";
"x1";"x2" **\dir{-};
"x2";"x3" **\dir{-};
"x3";"x4" **\dir{-};
"x4";"x1" **\dir{-};
(0,10)*{\bullet}="d2";
(0,20)*{\bullet}="d1";
(-2.6,10)*{d_2};
(-2.6,20)*{d_1};
(15,5)*{}="y1";
(25,5)*{}="y4";
(15,25)*{}="y2";
(25,25)*{}="y3";
"y1";"y2" **\dir{-};
"y2";"y3" **\dir{-};
"y3";"y4" **\dir{-};
"y4";"y1" **\dir{-};
(20,10)*{\bullet}="a2";
(20,20)*{\bullet}="a1";
(20,12.3)*{a_2};
(20,17.5)*{a_1};
(35,5)*{}="z1";
(45,5)*{}="z4";
(35,25)*{}="z2";
(45,25)*{}="z3";
"z1";"z2" **\dir{-};
"z2";"z3" **\dir{-};
"z3";"z4" **\dir{-};
"z4";"z1" **\dir{-};
(40,10)*{\bullet}="c2";
(40,20)*{\bullet}="c1";
(40.5,12.5)*{c_2};
(40.5,17.5)*{c_1};
(55,5)*{}="w1";
(65,5)*{}="w4";
(55,25)*{}="w2";
(65,25)*{}="w3";
"w1";"w2" **\dir{-};
"w2";"w3" **\dir{-};
"w3";"w4" **\dir{-};
"w4";"w1" **\dir{-};
(60,10)*{\bullet}="b2";
(60,20)*{\bullet}="b1";
(60.5,12.3)*{b_2};
(60.5,17.5)*{b_1};
"d2";"b2" **\dir{-};
"d1";"b1" **\dir{-};
"d1";"a2" **\dir{-};
"a1";"c2" **\dir{-};
"c1";"b2" **\dir{-};
"d1";"b1" **\crv{(30,32)} ?>* \dir{-};
"d2";"b2" **\crv{(30,-2)} ?>* \dir{-};
"d2";"c1" **\crv{(-40,45)} ?>* \dir{-};
"d2";"b1" **\crv{(0,-10) & (80,-20) & (80,10)} ?>* \dir{-};
"b2";"a1" **\crv{(100,45)} ?>* \dir{-};
"d1";"c2" **\crv{(-40,-15)} ?>* \dir{-};
"a2";"b1" **\crv{(100,-15)} ?>* \dir{-};
(20,-2)*{}="g4";
(67.5,-2)*{}="g5";
(67.5,20)*{}="g6";
\endxy
\]
\caption{The final configuration}\label{fi:4ed11}
\end{figure}

Now we have the CME$(a_1-c_2, c_1)$, unless $c_2-d_1$. This implies $c_1-d_2$.  Further, we get CME$(a_1-c_2,a_2)$ unless $a_1-b_2$, which implies that $a_2-b_1$. This accounts for all edges (see Figure \ref{fi:4ed11}).

 But now we have the CME$(b_2-d_2, b_1)$, as $N(b_2) \cap N(d_2) \cap K =\{c_1\} \subseteq N(b_1)$, and $N_B' \subseteq N(b_1)$. This contradiction
 excludes the final case in the proof of Theorem \ref{th:n-4}.

 We have shown that every primitive group synchronizes every transformation of rank $n-4$.

\section{Primitive groups of permutation rank~$3$}

The arguments of the preceding sections apply in complete generality because
they use only the fact that the groups involved are primitive. We can get stronger
results by focussing on a restricted class of primitive groups, in particular, the
primitive permutation groups of rank~$3$. (Unfortunately the term
``rank'' is used in a different sense by permutation group theorists!)

More precisely, the (permutation) \emph{rank} of a transitive permutation group
$G$ acting on a set $X$ is the number of orbits of $G$ on $X\times X$, the set of ordered
pairs of elements of $X$. Equivalently, it is the number of orbits on $X$ of
the stabiliser of a point of $X$.

If $|X| > 1$, then the rank of $G$ is at least $2$, because
no permutation can map $(x,x)$ to $(x,y)$. A primitive group of rank~$2$
is doubly transitive (and hence synchronizing), and thus the first non-trivial cases
are primitive groups of rank~$3$. The aim of this section is to prove the following result. 

\begin{theorem}\label{Artur} A primitive permutation group of degree $n$ and permutation  rank~$3$
synchronizes any map with rank strictly  larger than $n - (1+\sqrt{n-1}/12)$. 
\end{theorem}

Although a complete classification of the primitive  groups of rank~$3$ is known (see \cite{kl,ls,liebeck}),
we do not use this, but use instead combinatorial properties of strongly regular graphs. (A graph is
\emph{strongly regular} if the numbers $k$, $\lambda$, $\mu$ of neighbours of a vertex, an edge, and a non-edge respectively are independent of the chosen vertex, edge or non-edge. See~\cite{CvL} for the definition
and properties of strongly regular graphs. It is well known that a group with
permutation rank $3$ is contained in the automorphism group of a strongly regular graph.)

More precisely, we shall prove
the following result for strongly regular graphs. In the statement of this result -- and throughout this section -- we call a strongly regular graph {\em non-trivial} if it is connected and its complement is connected, which is the same as requiring that $\mu >  0$ and $k > \mu$ (the word ``primitive'' is sometimes
used to denote this property, but to avoid confusion with our many other uses of primitive, we will not use it in this sense).

\begin{theorem}\label{srg}
Let $\Gamma$ be a non-trivial strongly regular graph on $n$ vertices and let $f \in \mathrm{End}(\Gamma)$
be an endomorphism of $\Gamma$ of rank $r$. Then $n - r \geq 1 + \sqrt{n-1}/12$.
\end{theorem}

The proof of this uses three simple lemmas:
\begin{lemma}\label{n-rbound}
If $\Gamma$ is a non-trivial strongly regular graph with parameters $(n,k,\lambda,\mu)$, and $f$ is a proper endomorphism of $\Gamma$ of
rank $r$, then
\begin{equation*}
n-r \geq (k-\mu+4)/4.
\end{equation*}
\end{lemma}

\begin{pf}
Suppose that the kernel of $f$ has $t$ singleton classes, and therefore $n-t$ vertices in
non-singleton classes. As $f$ is not an automorphism, it follows that $n-t \geq 2$, and
because the non-singleton classes each have size at least $2$, we have $r \leq t + (n-t)/2$. By adding $(n-t)/2$ to each side
of this last expression and rearranging, we conclude that $n-t \leq 2(n-r)$.

Let $v$ and $w$ be two vertices in the same kernel class of $f$ and let $V$, $W$
be the neighbours of $v$ and $w$ respectively that lie in singleton kernel classes.
As $f$ identifies $v$ and $w$, and maps the vertices of $V \cup W$ injectively to the neighbours of $vf$
it follows that $|V \cup W| \leq k$.
Vertices $v$ and $w$ are each adjacent to at most $n-t-2$ vertices lying in non-singleton kernel classes so
$|V| \geq k-(n-t-2)$ and similarly for $|W|$.
Therefore
\begin{align*}
|V \cap W| &= |V| + |W| - |V \cup W| \\
&\geq k-(n-t-2) + k-(n-t-2) - k \\
&= k - 2(n-t) + 4\\
& \geq k - 4(n-r) + 4.
\end{align*}
Finally, as $v$ and $w$ are not adjacent, it follows that $|V \cap W| \leq \mu$ and the result follows
by combining the two bounds for $|V \cap W|$. \qed
\end{pf}

\begin{lemma}\label{k-mubound}
If $\Gamma$ is a non-trivial strongly regular graph with parameters $(n,k,\lambda,\mu)$, then
\begin{equation*}
k-\mu\geq \frac{1}{3}\min(k,k').
\end{equation*}
where $k' = n-k-1$ is the valency of the complement of $\Gamma$.
\end{lemma}
\begin{pf}
If $\Gamma$ is a conference graph, then $n=4\mu+1$ and $k=2\mu$ and so $k-\mu= k/2 = k'/2$, thereby satisfying the
conclusion of the theorem. Otherwise the three eigenvalues of $\Gamma$, which we denote $k$, $r$ and $s$ (with $r > 0 > s$), are
all integers, and in particular $r \geq 1$. (There is possible confusion with
the use of $r$ as the rank of an endomorphism; note that we only use $r$ in
the present sense within this proof, following the notation of \cite{CvL},
and endomorphisms will not occur here.)

It is well-known that all the parameters of a strongly regular graph can be expressed purely in  terms
of $k$, $r$ and $s$ (see~\cite[Chapter 2]{CvL}) and from this it can be
deduced that
\[\frac{kr(k'+r+1)}{k(r+1)+k'r}=\frac{krs(r+1)(r-k)}{k(k-r)(r+1)}=-rs,\]
by substituting
\[k'=\frac{k(k-\lambda-1)}{\mu}=\frac{-k(r+1)(s+1)}{k+rs}\]
into the left-hand side.
From this, we can conclude that
\begin{align*}
k-\mu &= -rs = \dfrac{k(k'+r+1)}{k(1+\frac{1}{r})+k'} \geq \begin{cases}
\dfrac{k'}{2+\frac{k'}{k}}\geq \dfrac{1}{3} k', \text{ for } k'\leq k.\\
\\
\dfrac{k}{2\frac{k}{k'}+1} \geq  \dfrac{1}{3}k, \text{ for } k\leq k'.
\end{cases}
\end{align*}
where the final inequalities arise from dividing by either $k$ or $k'$, and then using the fact that $r \geq 1$. \qed
\end{pf}

\begin{lemma}\label{minkbound}
If $\Gamma$ is a non-trivial strongly regular graph with parameters $(n,k,\lambda,\mu)$, then
\begin{equation*}
\min(k, k') \geq \sqrt{n-1}.
\end{equation*}
\end{lemma}

\begin{pf}
As $\Gamma$ and its complement are both connected graphs of diameter~$2$,
the Moore bound implies that $n\le k^2+1$ and $n\le k'^2+1$ and the result
follows immediately. \qed
\end{pf}

Thus combining the results of Lemmas~\ref{n-rbound}, \ref{k-mubound} and \ref{minkbound}, we conclude that a
proper endomorphism of rank $r$ of a non-trivial strongly regular graph on $n$ vertices
satisfies
\[
n-r \geq 1 + \sqrt{n-1}/12,
\]
thereby completing the proof of Theorem~\ref{srg}.

\paragraph{Remark} The constant $1/12$ in this theorem is not best possible,
and can be improved by using the classification of primitive permutation groups of rank $3$ mentioned above.
Details will appear elsewhere.

\paragraph{Remark} No non-trivial strongly regular graphs are known
that have any proper endomorphisms other than colourings (i.e. endomorphisms whose
image is a clique).

\section{Computational Results}\label{compute}

In this section we briefly describe the results of searching for endomorphisms in small vertex-primitive graphs,
namely those on (strictly) fewer than $45$ vertices. In addition to confirming that the linegraph of the
Tutte-Coxeter graph is the smallest example of a vertex-primitive graph admitting a non-uniform endomorphism, there
are various points in the theoretical arguments that terminate by requiring that certain small cases
be checked, so for convenience, we gather all this information in one place.

The \emph{primitive groups} of small degree are easily available in both \gap\ and \magma, though the
reader is warned that these two computer algebra systems use \emph{different numbering systems} so
that, for example, \texttt{PrimitiveGroup(45,1)} is $\mathrm{PGL}(2,9)$ in \gap, but $M_{10}$ in
\magma. As we are only seeking vertex-primitive graphs whose chromatic number and clique number are
equal, we need not consider the primitive groups of prime degree, which have a large number of
orbitals. The remaining groups have a much more modest number of orbitals and it is easy to
construct all possible graphs stabilised by each group by taking every possible subset of the orbitals
(ensuring that if a orbital that is not self-paired is chosen, then so is its partner).

For the sizes we are considering (up to $45$ vertices), it is fairly easy to determine the
chromatic and clique numbers of the graphs and thus extract all possible graphs whose
endomorphism monoids might contain non-uniform endomorphisms. There
are only $24$ such graphs on fewer than $45$ vertices and in Table~\ref{omchi}, we give summary
data listing just the order $n$, the valency $k$ and the chromatic number $\chi$ of each of these graphs. For example, the entry $(12,5)^3$
in the row for $n=25$ indicates that on $25$ vertices, there are three $12$-regular vertex-primitive graphs with $\omega=\chi = 5$.
There are no further examples on $37$--$44$ vertices and so this list is complete for $n < 45$.

\begin{table}
\begin{center}
\begin{tabular}{cl}
$n$ & \multicolumn{1}{c}{Values of $(k,\chi)$ occurring}\\
\hline
$9$ & $(4,3)$\\
$15$& $(8,5)$\\
$16$& $(6,4)$, $(9,4)$\\
$21$& $(4,3)$, $(16,7)$\\
$25$& $(8,5)$, $(12,5)^3$, $(16,5)$\\
$27$& $(6,3)$, $(8,3)$, $(18,9)$, $(20,9)$\\
$28$& $(6,4)$, $(12,7)$, $(15,7)$, $(18,7)^2$, $(21,7)$\\
$35$& $(18,7)$\\
$36$& $(10,6)$, $(25,6)$\\
\hline
\end{tabular}
\end{center}
\caption{$(k,\chi)$ for $n$-vertex primitive graphs with $\omega=\chi$}\label{omchi}
\end{table}

The bottleneck in this process is not the construction of the graphs, nor the calculation
of their chromatic or clique numbers, but rather the computation of their
endomorphisms. Apart from some obvious use of symmetry (for example, requiring that
a vertex be fixed), we know no substantially better method than to perform what is essentially a
naive back-track search. This finds an endomorphism by assigning to each vertex in turn
a candidate image, determines the consequences of that choice (in terms of reducing the possible choices
for the images of other vertices), and then turns to the next
vertex, until either a full endomorphism is found, or there are unmapped vertices for which
no possible choice of image respects the property that edges are mapped to edges.

Such a search can easily be programmed from scratch, but in this case we used the
constraint satisfaction problem solver {\sc Minion}. This software, which was developed at St Andrews,
performs extremely well for certain types of search problem. Using {\sc Minion}, we confirmed
that for all but two of the graphs listed in Table~\ref{omchi}, every endomorphism is
either an automorphism or a colouring. The two exceptions are the $6$- and $8$-regular graphs on
$27$ vertices which also have ``in-between'' endomorphisms whose image is the $9$-vertex Paley graph $P(9)$.
The $6$-regular graph is the Cartesian product $P(9) \cart K_3=K_3\cart K_3\cart K_3$,
while the $8$-regular graph is the direct product $P(9) \times K_3=K_3\times K_3\times K_3$.

On $45$ vertices, there are $8$ non-trivial vertex-primitive graphs with equal chromatic and clique number,
including the linegraph of the Tutte-Coxeter graph.  Of the remaining graphs, some are
sufficiently dense that we have been unable yet to completely determine all of their endomorphisms.
However by a combination of computation and theory, we at least know that none of the $45$-vertex
graphs other than the linegraph of the Tutte-Coxeter graph admit proper endomorphisms other than colourings.

\section{Problems}\label{spro}
\setcounter{theorem}{0}

This paper started with the intention of providing further evidence that primitive groups are
almost synchronizing but, rather inconveniently, this turns out not to be true. Therefore, faced
with an unexpectedly complex situation, we pose the following problem, although with the expectation
that resolving it is likely to be difficult:

\begin{prob}
Classify the almost synchronizing primitive groups.
\end{prob}

It might be more feasible to focus on the ``large-rank'' end of the spectrum,
where we still believe that the following weaker version of the almost synchronizing conjecture
is true.

\begin{conjecture}\label{c:half}
A primitive group of degree $n$ synchronizes any map whose rank $r$
satisfies $n/2<r<n$ (all such maps are non-uniform).
\end{conjecture}

As we have seen, showing that primitive groups synchronize maps of rank
$n-4$ required a long case analysis. Further progress will require a solution
of the following problem.

\begin{prob}
Find new techniques to show that large-rank transformations are synchronized
by primitive groups, and use them to extend the range below $n-4$.
\end{prob}

The previous problems deal with the spectrum of ranks synchronized by primitive groups. An orthogonal approach is to investigate the kernel types that are synchronized by primitive groups, along the line of the results in Section \ref{k-2sets}.

\begin{prob}
Find new kernel types synchronised by a primitive group. In particular, prove that all primitive groups synchronize maps with the following kernel types: 
\[	
(2,\ldots,2,1,\ldots,1) \mbox{ or  } (p,q,1,\ldots,1), \mbox{ for all $p,q>1$ .}
\]	
 \end{prob}
 
\begin{prob}
Is there a ``threshold'' function $f$ such that a transitive permutation group
of degree $n$ is imprimitive if and only if it has more than $f(n)$
non-syn\-chronizing ranks? (A positive answer to Conjecture~\ref{c:half} would
show that $f(n)=n/2$ would suffice.) In particular, is the number of
non-synchronizing ranks of a primitive group $o(n)$?
\end{prob}

Theorem \ref{Artur} concerns the synchronizing power  of groups of permutation rank $3$ and maps of \emph{large rank}. In the spirit of the remaining results of this paper, it would be interesting to investigate what happens with maps of \emph{small rank}. 
\begin{prob}
Find the largest natural number $k$ such that groups of permutation rank $3$ synchronize every non-uniform map of rank $l$, for all $l\le k$. 
\end{prob}

The previous problem is somehow connected to the next, the classification of a class of groups  lying strictly between primitive and synchronizing.

\begin{prob}
Is it possible to classify the primitive groups which synchronize
every rank $3$  map?
\end{prob}
	
	The previous problem is equivalent to classifying the permutation groups $G$, acting primitively on a set $\Omega$, such that for every $3$-partition $P$ of $\Omega$ and every section $S$ for $P$, there exists $g\in G$ such that $Sg$ is not a section for $P$. 	
	
Note that there are primitive groups that do not synchronize a rank $3$ map (see the example immediately before Section \ref{trans} in \cite{arcameron22}). And there are non-synchronizing groups which synchronize every rank $3$ map. Take for example $\pgl(2,7)$ of degree $28$; this group is non-synchronizing, but synchronizes every rank $3$  map	since $28$ is not divisible
by $3$.

There are very fast algorithms to decide if a given set of permutations generate a primitive group, but is it possible that such
an algorithm exists for synchronization?

\begin{prob}
Find an efficient algorithm to decide if a given set of permutations generates a synchronizing group or show that
such an algorithm is unlikely to exist.
\end{prob}

\begin{prob}\label{11}
Formulate and prove analogues of our results for semigroups of linear maps on a
vector space. Note that linear maps cannot be non-uniform, but we could ask for
linear analogues of results expressed in terms of rank such as Theorem \ref{main}.
\end{prob}

\begin{prob}
Solve the analogue of Problem \ref{11} for independence algebras (for definitions and fundamental results see 
\cite{Ar1,ArEdGi,arfo,Ar3,abk,cameronSz,F1,F2,gould})
\end{prob}

 Suppose the diameter of a group $G$ (acting on a set $\Omega$) is at most $n-1$ (that is, given any set $S$ of generators of $G$, every element of $G$ can be generated by the elements of $S$ in a word of length at most $n$). Suppose, in addition, that $G$ and a transformation $t$ of $\Omega$ generate a constant map $tg_1t\ldots g_{n-2}t$. Then we can replace the $g_i$ by a word (on the elements of $S$) of length at most $n$ and hence we have a constant written as a word of length meeting the \v{C}erny bound. However, finding the diameters of primitive groups is a very demanding problem. Therefore we suggest the following two problems.

\begin{prob}
Let $\Omega$ be a set. Let $G$ be a synchronizing group acting primitively on $\Omega$ and let $S\subseteq G$ be a set of generators for $G$. Let $X\subseteq \Omega$ be a proper subset of $\Omega$, and let $P$ be a partition of $\Omega$ in $|X|$ parts. Is it true that there exist two elements in the set $X$ that can be carried to the same part of $P$ by a word (on the elements of $S$) of length at most $n$?
\end{prob}

 We consider the previous problem one of the most important by its implications on  the \v{C}erny conjecture, in the case of transformation semigroups containing a primitive synchronizing group. 
 
 The previous problem admits also a general version for primitive groups.

\begin{prob}
Let $\Omega$ be a set. Let $G$ be a group acting primitively on $\Omega$ and let $S\subseteq G$ be a set of generators for $G$. Let $X\subseteq \Omega$ be a proper subset of $\Omega$, and let $P$ be a partition of $\Omega$ in $|X|$ parts. Let $Q\subseteq X\times X$ be the set of pairs $(x,y)$ such that for some $g\in G$ we have $xg$ and $yg$ belonging to the same part of $P$. Assuming $Q\neq \emptyset$, is it true that there exists $(x_0,y_0)\in Q$ and a word $w$ (on the elements of $S$), of length at most $n$, such that $x_0w$ and $y_0w$ belong to the same part of $P$?
\end{prob}

The computations in this paper were critical to prove our results; and the generalizations of our results will certainly require to push the limits of the computations above.

\begin{prob}
Extend the computational results of Section \ref{compute}.
\end{prob}
\section*{Acknowledgements}
The second author  has received funding from the
European Union Seventh Framework Programme (FP7/2007-2013) under
grant agreement no.\ PCOFUND-GA-2009-246542 and from the Foundation for
Science and Technology of Portugal under  PCOFUND-GA-2009-246542 and SFRH/BCC/52684/2014, as well as through the CAUL / CEMAT project.

The first author has been partially supported by an grant from the Foundation for
Science and Technology of Portugal.

\end{document}